\pdfoutput=1
\RequirePackage{ifpdf}
\ifpdf 
\documentclass[pdftex]{sigma}
\else
\documentclass{sigma}
\fi

\numberwithin{equation}{section}

\newtheorem*{Theorem*}{Theorem}

 { \theoremstyle{definition}

 }

\usepackage{scalerel}
\newcommand\bigtimes{\mathbin{\ThisStyle{\vcenter{\hbox{%
\scalebox{1.5}{$\SavedStyle\times$}}}}}}
\begin{document}

\allowdisplaybreaks

\renewcommand{\thefootnote}{}

\newcommand{\arXivNumber}{2306.02932}

\renewcommand{\PaperNumber}{038}

\FirstPageHeading

\ShortArticleName{Product Inequalities for $\mathbb T^\rtimes$-Stabilized Scalar Curvature}

\ArticleName{Product Inequalities for $\boldsymbol{\mathbb T^\rtimes}$-Stabilized Scalar\\ Curvature\footnote{This paper is a~contribution to the Special Issue on Differential Geometry Inspired by Mathematical Physics in honor of Jean-Pierre Bourguignon for his 75th birthday. The~full collection is available at \href{https://www.emis.de/journals/SIGMA/Bourguignon.html}{https://www.emis.de/journals/SIGMA/Bourguignon.html}}}

\Author{Misha GROMOV~$^{\rm ab}$}

\AuthorNameForHeading{M.~Gromov}

\Address{$^{\rm a)}$~Courant Institute of Mathematical Sciences, New York University,\\
\hphantom{$^{\rm a)}$}~New York, NY~10012-1185, USA}

\Address{$^{\rm b)}$~Institut des Hautes \'Etudes Scientifiques, 91893 Bures-sur-Yvette, France}
\EmailD{\href{mailto:gromov@ihes.fr}{gromov@ihes.fr}}
\URLaddressD{\url{https://www.ihes.fr/~gromov/}}

\ArticleDates{Received June 26, 2023, in final form April 26, 2024; Published online May 08, 2024}

\Abstract{We study metric invariants of Riemannian manifolds $X$ defined via the $\mathbb T^\rtimes$-{\it sta\-bi\-lized scalar curvatures} of manifolds $Y$ mapped to $X$ and prove in some cases additivity of these invariants under Riemannian products $X_1\times X_2$.}


\Classification{53C21}

\begin{flushright}
\begin{minipage}{48mm}
\it To the 75th birthday\\ of Jean-Pierre Bourguignon
\end{minipage}
\end{flushright}

\renewcommand{\thefootnote}{\arabic{footnote}}
\setcounter{footnote}{0}

\section[T\^{}\{rtimes\}-stabilization]{$\boldsymbol{\mathbb T^\rtimes}$-stabilization}

A {\it ``warped'' $\mathbb T^N$-extension}, $N=0,1,\dots$,
of a Riemannian manifold $X=(X,g)$, possibly with a~boundary, is
\[
X_N^\rtimes=X\rtimes \mathbb T^N = \bigl(X\times \mathbb T^N, g^\rtimes\bigr),
\]
 where $\mathbb T^N$ is the (flat split) $N$-torus and where $X_N^\rtimes$ is endowed with a {\it warped metric},
\[
g^\rtimes= g_N^\rtimes =g_{N,\{\varphi_i\}}^\rtimes=g+\sum_{i=1}^N\varphi_i^2\,{\rm d}t_i^2,
\]
where
 $\varphi_i(x)\geq 0$ are smooth positive functions, which are strictly positive $(>0)$ in the interior $X\setminus \partial X$ of $X$.

Assume $g$ is smooth and let $\operatorname{Sc}^\rtimes_{\{\varphi_i\}} (X)$ be the scalar curvature of $g^\rtimes$, that is,
\[
\operatorname{Sc}^\rtimes_{\{\varphi_i\}} (X)=\operatorname{Sc}(g_{N,\{\varphi_i\}}^\rtimes)=\operatorname{Sc}\left(g+\sum_{i=1}^N\varphi_i^2\,{\rm d}t_i^2\right),
 \]
 where $\operatorname{Sc}^\rtimes_{\{\varphi_i\}} (X)=\operatorname{Sc}^\rtimes_{\{\varphi_i\}} (X,x)$ is a function on $X$, since $g_N^\rtimes$ is invariant under the obvious action of $\mathbb T^N$ on $X_N^\rtimes=X\times \mathbb T^N$.

Let $\operatorname{Sc}^\rtimes (X)$, $X=(X,g)$, be the supremum of the numbers $\sigma$ such that
 $\operatorname{Sc}_{N,\{\varphi_i\} }^\rtimes (X)>\sigma$ for some $N$ and $ \varphi_i$.

\textbf{$\boldsymbol{\supset}$-Monotonicity.} Clearly,
\[
\operatorname{Sc}^\rtimes (Y)\geq \operatorname{Sc}^\rtimes (X)
\]
 for all smooth domains $Y\subset X$.

\textbf {1.A. Formulas.} A computation shows\footnote{See \cite{F-CS1980}, \cite[formulas~(7.33) and (12.5)]{GL1983}, also \cite{Zhu2019}, and  \cite[Section~2.4.1]{Gr2021}.} that
\[
\operatorname{Sc}\left(g+\sum_{i=1}^N\varphi_i^2\,{\rm d}t_i^2\right)= \operatorname{Sc}(g)-2\sum^N_{i=1} \frac{\Delta\varphi_i} {\varphi_i}-2\sum_{i<j}
 \langle\nabla_g\log\varphi_i,\nabla_g\log\varphi_j
 \rangle.
 \]
 For instance, if $N=1$ and $\varphi_1$ is a positive eigenfunction of the operator
 $-\Delta+\frac{1}{2}\sigma$, for
\[
\sigma(x)=\operatorname{Sc}(g,x),
\]
that is,
\[
\varphi(x)\mapsto -\Delta \varphi(x) +\frac{1}{2}\sigma(x)\varphi(x)
\]
 and
\[
-\Delta \varphi_1 +\frac{1}{2}\sigma(x)\varphi_1=\lambda_1\varphi_1,
\]
 then
 \[
 \operatorname{Sc}\bigl(g+\varphi_1^2\,{\rm d}t^2\bigr)
 =\sigma -2\frac {\Delta\varphi_1} {\varphi_1}=\sigma -2 \frac{\frac{1}{2}\sigma\varphi_1- \lambda \varphi_1}{
 \varphi_1} =2\lambda_1.
 \]
Thus,
\[
\operatorname{Sc}^\rtimes (X)\geq 2\lambda_1=2\lambda_1\left(\Delta+\frac{1}{2}\operatorname{Sc}(x)\right).
\]

Recall at this point that if $X$ is a compact connected manifold without boundary, then, for all functions $\sigma(x)$,
\begin{align*}
\lambda_1& =\inf_{\varphi\neq 0}\left(\frac {-\int_X\varphi(x) \Delta \varphi(x){\rm d}y} {\int_X\varphi(x)^2{\rm d}y}+
 \frac12\frac {\int_X
 \sigma(x)\varphi^2(x){\rm d}y}{\int_X\varphi(x)^2{\rm d}y} \right)\\
 & \geq
\inf_{\varphi\neq 0}\left(\frac {-\int_X\varphi(x) \Delta \varphi(x){\rm d}y} {\int_X\varphi(x)^2{\rm d}y}+\frac12
 \inf_{x\in X}\sigma(x)\right) \geq \frac12\inf_{x\in X}\sigma(x),
\end{align*}
 where the latter inequality follows from positivity of $-\Delta$ and this inequality is strict $(>)$ with $\varphi=\varphi_1$, unless~$\sigma(x)$ is constant.

Also, the strict inequality $\lambda_1> \inf_{x\in X}\operatorname{Sc}(X,x)$ holds for the {\it Dirichlet eigenvalue} $\lambda_1$
 on compact connected manifolds {\it with boundaries},
 since the above relations
 are satisfied for functions~$\varphi (x)$, which vanish on the boundary.

 Next, now for all $N$ and all $\varphi_i$, rewrite the above expression for $\operatorname{Sc}\bigl(g+\sum_{i=1}^N\varphi_i^2{\rm d}t_i^2\bigr)$
 with the function $\Phi(x)=\log(\varphi_1(x)\cdots \varphi(x)_N)$ as follows:
\[
\operatorname{Sc}\left(g+\sum_{i=1}^N\varphi_i^2\,{\rm d}t_i^2\right)= \operatorname{Sc}(g) -2\Delta\Phi-||\nabla \Phi||^2 -\sum_i{||\nabla\log \varphi_i||^2}.
\]
This shows that
$\operatorname{Sc}\bigl(g+\sum_{i=1}^N\varphi_i^2\,{\rm d}t_i^2\bigr)$ {\it increases} under replacing all $\varphi_i$ by their {\it geometric mean},
$\varphi_i\leadsto \phi= \sqrt[N]{\prod_i\varphi_i}$, i.e.,
\[
\operatorname{Sc}(g)-2\Delta\Psi -\frac {N+1}{N}||\nabla\Psi||^2\geq \operatorname{Sc}(g) -2\Delta\Phi-||\nabla \Phi||^2 -\sum_i||\nabla\log\varphi_i||^2
\]
for $\Psi=\log \phi^N$, where the equality holds only if all $\nabla\log\varphi_i$ are mutually equal.
Hence,
\[
\sup_{\varphi_i>0} \operatorname{Sc}\left(g+\sum_{i=1}^N\varphi_i^2\,{\rm d}t_i^2\right)=\sup_\Psi \operatorname{Sc}(g)-2\Delta\Psi -\frac {N+1}{N}||\nabla\Psi||^2,
\]
where this ``$\sup$'' {\it increases} with $N$; thus,
 by letting $N\to\infty$, we see that
\begin{align*}
\operatorname{Sc}^\rtimes (X)& =\sup_{\Psi(x) }\inf_{x\in X} \operatorname{Sc}(X,x)-2\Delta\Psi(x) -||\nabla\Psi(x)||^2 \\
& = \sup_{\psi(x)>0}\inf_{x\in X} \operatorname{Sc}(X,x)- 2\frac{\Delta\psi(x)} {\psi(x)} +\frac{||\nabla\psi(x)||^2}{\psi(x)^2}.
\end{align*}
Rewrite this equation with $\Psi=2\Theta$ as follows:
\begin{align*}
\operatorname{Sc}^\rtimes (X) & =\sup_{\Theta }\inf_{x\in X} \operatorname{Sc}(X,x)-4\bigl(\Delta\Theta(x) +||\nabla\Theta(x)||^2\bigr)\\
& =
\sup_{\theta }\inf_{x\in X} \operatorname{Sc}(X,x)-4\frac{\Delta\theta(x)} {\theta(x)} \qquad \text{for} \quad \theta=\exp \Theta.
\end{align*}

 Therefore, if $X$ is compact, then
\[
\operatorname{Sc}^\rtimes (X)\geq 4\lambda_1^\rtimes(X),
\]
where $\lambda^\rtimes_1(X)$ is the lowest eigenvalue of the operator
$-\Delta +\frac {1}{4}\operatorname{Sc}$ on $X$ with the Dirichlet (vanishing on the boundary) condition.

(If a connected manifold $X$
has no boundary and the scalar curvature of $X$ is constant,
then $\operatorname{Sc}^\rtimes (X)=\operatorname{Sc}(X)$; otherwise
\[
\operatorname{Sc}^\rtimes (X)> \inf_{x\in X}\operatorname{Sc}(X,x).
\]
 Moreover,
\[
\operatorname{Sc}^\rtimes (X)>\beta^{-1} \lambda_1\left(-\Delta +\beta \operatorname{Sc}(X)\right)
\]
 for all $\beta>1/4$, since the operator $-\Delta$ is
{\it strictly positive} on non constant functions on $X$.)

\textbf{1.B. $\boldsymbol{\frac {1}{4}}$-Proposition.} Let $X=(X,g)$ be a compact Riemannian manifold with a boundary. Then
\[
\operatorname{Sc}^\rtimes (X)= 4\lambda_1^\rtimes(X).
\]
In fact, let $\theta_1(x)=\theta^\rtimes_1(x)\geq 0$, be the first Dirichlet eigenfunction of the operator $-\Delta +\frac {1}{4}\operatorname{Sc}$ and
let $\theta(x)>0$ be an arbitrary smooth strictly positive function on $X$.

Since $\theta_1$ is strictly positive in the interior of $X$,
the ratio \smash{$\frac{\theta_1(x)}{\theta(x)}$} assumes its maximum, call it~$a$, at an interior point $x_0\in X$, where
\[
\theta_1(x_0)=a\theta(x), \qquad a \nabla\theta_1(x_0)= \nabla\theta(x) \qquad \text{and} \qquad \Delta\theta_1(x)\leq a\Delta\theta(x),
\]
and consequently
\[
-\frac{\Delta\theta_1(x_0)} {\theta_1(x_0)}+\frac14 \operatorname{Sc}(X,x_0)\geq -\frac{\Delta\theta(x_0)} {\theta(x_0)}+\frac14 \operatorname{Sc}(X,x_0).
\]
Since the sum \smash{$-\frac{\Delta\theta_1(x)} {\theta_1(x)}+\frac14 \operatorname{Sc}(X,x)$} is constant ($=\lambda_1^\rtimes$), this
inequality holds at the minimum point $x_o\in X$ of the function
 \smash{$ -\frac{\Delta\theta(x)} {\theta(x)}+\frac14 \operatorname{Sc}(X,x)$}; hence,
\smash{$\inf_{x\in X} \operatorname{Sc}(X,x)-4\frac{\Delta\theta\rtimes_1(x)} {\theta^\rtimes (x)}$} {\it majorizes} \smash{$\inf_{x\in X} \operatorname{Sc}(X,x)-4\frac{\Delta\theta_1(x)} {\theta (x)}$} for all smooth {\it strictly positive} functions on~$X$.

Finally, by {\large$\supset$}-monotonicity of $\operatorname{Sc}^\rtimes$ and the continuity of the first eigenvalue $\lambda_1(-\Delta_g + \operatorname{Sc}(g)/4)$ in (the space of $C^2$-metrics) $g$, this majorization holds for functions $\theta$ which is strictly positive {\it only} in the interior of $X$.\footnote {One needs to be slightly careful here, since
$\Delta \theta/\theta$ may, a priori, blow up at the boundary of $X$.} Then
the proof of the $\frac {1}{4}$-proposition follows.

\textbf{$\boldsymbol{\frac{1}{4}}$-Remark.} This $\frac{1}{4}$ agrees with that in the
Schr\"odinger--Lichnerowicz formula $\mathcal D^2=\nabla^\ast\nabla+\frac{1}{4}\operatorname{Sc}$ via the Kato inequality for the squared Dirac operator on $X^\rtimes =\bigl(X\times \mathbb T^N, g^\rtimes\bigr)$, which, to make it index-wise more interesting, may be twisted with the canonical $N$-parametric family of flat unitary complex line bundles over~$X^\rtimes$.
(Probably, there is much of what we do not understand about the relations between the two~$\frac{1}{4}$.)

\subsection*{Corollaries/Examples}

\textbf{1.B$\boldsymbol{{}_1}$.} If $X$ is {\it non-compact}, then
\[
\operatorname{Sc}^\rtimes(X) =\lim_{i\to\infty} \operatorname{Sc}^\rtimes(X_i)
\]
for compact equidimensional submanifolds
$X_1\subset \cdots \subset X_i\subset \cdots \subset X,$
which {\it exhaust} $X$.

 \textbf{1.B$\boldsymbol{{}_2}$.} $\operatorname{Sc}^\rtimes $ is additive under Riemannian products:
\[
\operatorname{Sc}^\rtimes (X_1\times X_2)=\operatorname{Sc}^\rtimes (X_1) + \operatorname{Sc}^\rtimes (X_2).
\]

For instance, the rectangular solids satisfy
\[
\operatorname{Sc}^{\rtimes}\left( \bigtimes_1^n[-a_i,b_i]\right)=4\sum_1^n \lambda_1[a_i,b_i]=
\sum_1^n \frac {4\pi^2}{(b_i-a_i)^2}.
\]

 \textbf{1.B$\boldsymbol{{}_3}$.} Manifolds $X$ with constant scalar curvature $\sigma$ satisfy
\[
\operatorname{Sc}^\rtimes(X)=4\lambda_1(X)+\sigma
\]
 for the first eigenvalue $\lambda_1$ of the Laplace operator on $X$.

 For instance,
 unit hemispheres satisfy
\[
\operatorname{Sc}^\rtimes  (S^n_+ )=n(n-1)+4n=n(n+3)
\]
 and the unit balls $B^n=B^n(1)\subset \mathbb R^n$ satisfy
\[
\operatorname{Sc}^\rtimes (B^n)=4j_\nu^2,
\]
 for the first zero of the Bessel function $J_\nu$, $\nu=\frac {n}{2}-1$, where
 $j_{-1/2}=\frac{\pi}{2}$, $j_0=2.4042\dots$, $j_{1/2}=\pi$ and if $\nu>1/2$, then{\samepage
\[
\nu+\frac{a \nu^\frac{1}{3}}{2^\frac{1}{3}}<j_\nu<
\nu+\frac{a \nu^\frac{1}{3}}{2^\frac{1}{3}}+\frac{3}{20}\frac{2^\frac{2}{3}a^2}{\nu^\frac{1}{2}},
\]
 where $a=\bigl(\frac {9\pi}{8} \bigr)^\frac{2}{3}(1+\varepsilon)\approx 2.32$ with
 $\varepsilon< 0.13\bigl(\frac {8}{2.847\pi}\bigr)^{2}$ \cite{QW1999}.}

Specifically,
\begin{gather*}\operatorname{Sc}^\rtimes \bigl(B^2\bigr)=4(2.404\dots)^2=23.116\ldots>10=\operatorname{Sc}^\rtimes \bigl(S^2_+\bigr),\\
\operatorname{Sc}^\rtimes \bigl(B^3\bigr)>36>18=\operatorname{Sc}^\rtimes \bigl(S^3_+\bigr),\\
\operatorname{Sc}^\rtimes \bigl(B^4\bigr)=4(3,817\dots)^2=52.727\ldots>28=\operatorname{Sc}^\rtimes \bigl(S^4_+\bigr),\\
\operatorname{Sc}^\rtimes \bigl(B^8\bigr)=4(6.380\dots)^2= 162.827\ldots>88 =\operatorname{Sc}^\rtimes \bigl(S^8_+\bigr).
\end{gather*}

 \textbf{1.B$\boldsymbol{{}_4}$. Ricci comparison inequality.} Let $X$ be a (metrically) complete Riemannian mani\-fold with a boundary
 such that $\operatorname{Ricci}(X)\geq (n-1)\kappa$
 and ${\rm mean.curv}(\partial x)\geq \mu$.
Then
\[
\operatorname{Sc}^{\rtimes}(X)\geq \operatorname{Sc}^\rtimes(\underline B^n_{\kappa, \mu}),
\]
where $\underline B^n_{\kappa, \mu}$ is the ball in the complete simply connected $n$-space
$\underline S^n_\kappa$ with sectional curvature
$\kappa$, and where the mean curvature of the boundary $\partial \underline B^n_{\kappa, \mu}$ is equal to $\mu$.\footnote{The corresponding, comparison inequality for the Dirichlet (Schr\"odinger) $\lambda_1(\Delta)$ (compare with \cite{CF1974} and \cite{BF2017})
has, undoubtedly, been known for at least 45 years and the relation
$\lambda_1=-\sup_{\varphi>0}\inf_{x\in X} \frac {\Delta\varphi(x)}{\varphi(x)}$ must be dated to
the 19th century. My apologies to the author(s), whose paper(s) I failed to find on the web.}

For instance,  if $\operatorname{Ricci}(X)\geq 0$ and $\mathrm{mean.curv}(\partial X)\geq n-1$, then
\[
\operatorname{Sc}^{\rtimes}(X)\geq \operatorname{Sc}^\rtimes( B^n=\underline B^n_{0,n-1})=4j_1^2.
\]

In fact, let $\varphi(b) =\phi(d(b))=\phi_{\kappa, \mu}(d)$ be the first Dirichlet eigenfunction in $B^n_{\rho, \mu}$ written as
a~function of $d=d(b)=\operatorname{dist}(b, \partial (\underline B^n_{\kappa, \mu}))$ and let
$\varphi(x)=\phi(\operatorname{dist}(X, \partial x))$. Then, by the Bishop comparison inequality,
\[
\frac {\Delta_X\varphi(x)}{\varphi(x)}
\geq \frac {\Delta_{\underline S^n_\kappa}\phi(d(b))}{\phi(d(b))}=\lambda_1(\Delta_{\underline B^n_\kappa}),
\qquad
 d(b)=\operatorname{dist}(X,\partial x),
\]
and the proof follows.

 \textbf{1.B$\boldsymbol{{}_5}$.} If $ B_{-1}^n(r)$ is the hyperbolic $r$-ball,
 then, clearly, $ \operatorname{Sc}^\rtimes (B_{-1}^n(r))$ is monotonically decreasing in $r$, asymptotically to
 $4\frac {j_{\nu}^2}{r^2}$ for $r\to 0$ and $ \operatorname{Sc}^\rtimes (B_{-1}^n(r))\to -(n-1)$ for $r\to \infty$.

 In fact, it follows from \cite[Theorem~3.3]{BF2017} that
\[
\operatorname{Sc}^\rtimes (B_{-1}^n(r))=-n(n-1)+(n-1)^2\left (\frac {1}{r^2} +c(r)\right )
\]
 for a bounded positive function $c(r)$ such that $c(r)\to 1$ for $r\to \infty$, and
\[
 \frac {1}{6}\leq c(r)\leq 1 \qquad \text{for $r\geq 1$ and $n\geq 2$}.
 \]
 Thus,
\[
\operatorname{Sc}^\rtimes (B_{-1}^n(r))>0 \qquad \text{for
 $r\leq \sqrt \frac{6(n-1)}{5n+1}$}, \qquad   \operatorname{Sc}^\rtimes (B_{-1}^n(r))<0 \qquad \text{for $r\geq 3$ and $n\geq 2$}
 \]
roughly.

\textbf{1.C.} \textbf{General torical stabilizations.} The most permissive torical ``extension'' of a Riemannian manifold $X$ is a Riemannian
 manifold $X^\natural $ with an isometric $\mathbb T^N$-action and an isometry $X^\natural/\mathbb T^N\leftrightarrow X$. Here,
 as earlier, one defines the number $\operatorname{Sc}^\natural(X)$, which is clearly $ \geq \operatorname{Sc}^\rtimes(X)$.

It seems, however -- I did not honestly checked this -- that the curvature formulas for Riemannian submersions \cite{O'N1966} imply that $\operatorname{Sc}^\natural(X)\leq \operatorname{Sc}^\rtimes(X)$.

 Alternatively, if
 the fibration $X^\natural \to X^\natural/\mathbb T^N=X$ {\it admits a section}, then the $\rtimes$-rendition of the Schoen--Yau argument\footnote{See~\cite{Gr2023} and references therein.} implies
the equality
\[
\operatorname{Sc}^\natural(X)= \operatorname{Sc}^\rtimes(X)
\]
for $\dim(X)=n\leq 8$,\footnote {The case $n=8$ depends on~\cite{Sm2003}.} while for all $n$
this may follow from~\cite{WXY2021}, where both arguments apply not only to Riemannian submersions but to all distance non-increasing maps $\bigl(X\times \mathbb T^N, G\bigr)\to (X,g)$.

\textbf{1.D. $\boldsymbol{\lambda_1(\beta)}$-Remarks.} The formal properties of the operators $-\Delta+ \beta \operatorname{Sc}$ are similar for all~$\beta$, e.g.,
if a Riemannian manifold $X_0$ admits a locally isometric equidimensional map to $X$, then
\[
\lambda_1(-\Delta_{X_0}+\beta \operatorname{Sc}(X_0)) \geq \lambda_1(-\Delta_{X}+\beta \operatorname{Sc}(X));
\]
the spectra $ \{\lambda_1,\lambda_2,\dots\}$ of the operators $-\Delta+ \beta \operatorname{Sc}$ are additive under Riemannian products,
\begin{gather*}
\operatorname{spec}(-\Delta_{X_1\times X_2} + \beta \operatorname{Sc}(X_1\times X_2))=\operatorname{spec}(-\Delta_{X_1} + \beta \operatorname{Sc}(X_1))+ \operatorname{spec}(-\Delta_{ X_2} + \beta \operatorname{Sc}( X_2)).
\end{gather*}

There are several special values of $\beta$:
\begin{itemize}\itemsep=0pt

\item if $\beta=(N+1)/4N$ $(>1/4)$, then the corresponding $\lambda_1$ is equal
to the maximal constant scalar curvature of the warped metrics on $X\times \mathbb T^N$;

\item if $ \beta=(n-1)/4n$ $(<1/4)$, $n=\dim(X)$, this implies positivity of the square of the Dirac operator by {\it refined Kato's inequality} \cite{He2000}.

\item if $\beta=(n-2)(4(n-1))$ $(<(n-1)/4n)$, then the inequality $\lambda_1>0$ implies the existence of a conformal metric on $X$ with $\operatorname{Sc}>0$ by the Kazdan--Warner theorem.
    \end{itemize}

The geometric meaning of other $\beta$, as well
as of the higher eigenvalues $\lambda_i(X,\beta) $ of $-\Delta + \beta \operatorname{Sc}$ is unclear.\footnote{The second eigenvalue of $\lambda_2(-\Delta_ + \frac{1}{2} \operatorname{Sc})$ is used by Marques and Neves in the proof of the $S^3$-min-max theorem,~\cite{MN2011}, but the role of $\lambda_i(X,\beta) $ remains problematic for $\dim(X)\geq 4$,
 $i\geq 2$ and all $\beta$.}

\section[Sc\^{}\{rtimes downarrow\}, Sc\_\{sp\}\^{}\{rtimes downarrow\}, ..., Sc\^{}\{rtimes downarrow\}\_\{ast\} on homology]{$\boldsymbol{\operatorname{Sc}^{\rtimes_\downarrow}, \operatorname{Sc}_{\rm sp}^{\rtimes_\downarrow}, \dots, \operatorname{Sc}^{\rtimes_\downarrow}_\ast}$ on homology}\label{section2}

 Let $X$ be a metric space, e.g., a Riemannin manifold $\operatorname{Sc}^{\rtimes_\downarrow}(h)=\operatorname{Sc}_{\rm dist}^{\rtimes_\downarrow}(h)$, $h\in H_m(X,\partial X)$, denote the supremum of the numbers $\sigma$ such that
the homology class $h$ is representable by
 a~{\it distance decreasing} map $f$ from an oriented Riemannian $m$-manifold $Y$ with $\operatorname{Sc}^\rtimes (Y) \geq \sigma$,
\[
f\colon \ (Y,\partial Y)\to (X, \partial X), \qquad
 f_\ast [Y,\partial Y]=h.
 \]

{\it Smoothness remark.} If $X$ is a smooth Riemannian manifold, then an obvious approximation argument shows that requiring maps $f$ to be smooth does not change the
value of $\operatorname{Sc}_{\rm dist}^{\rtimes_\downarrow}(h).$
(However, smoothness of a {\it distance non-increasing} map $f$ in the {\it extremal case}, where $\operatorname{Sc}^\rtimes (Y) =\operatorname{Sc}^{\rtimes_\downarrow}(h)$ is a delicate matter, see~\cite{CHS2023}.)

Below the are several versions of this definition with the generic notation $\operatorname{Sc}^{\rtimes_\downarrow}_\ast$.

\textbf{I.} Restrict/relax the topology of $Y$, e.g., by requiring that
\begin{description}\itemsep=0pt\setlength{\leftskip}{0.50cm}
\item[$\bullet_{\rm sp}$] $Y$ is spin;\footnote{A referee suggested to constrain maps
rather then manifolds $Y$, e.g., by allowing spin maps only (in the case, where $X$ is a manifold), but I could not figure out what to do with it.}
\item[$\bullet_{\widetilde{\rm sp}} $] the universal covering of $Y$ is spin;\footnote{This condition is satisfied in several interesting examples of {\it non-spin} manifolds $X$, e.g., where $\pi_2=0 $. At the same time, much of the Dirac theoretic scalar curvature results apply to these $X$, see \cite[Section~$9\frac12$]{Gr1996}.}
\item[$\bullet_{\pi_2=0}$] the second homotopy group of $Y$ is zero;
\item[$\bullet_{\rm st.par}$] $Y$ is stably parallelizable;
\item[$\bullet_{\odot}$] allow representation of $h$
 by {\it quasi-proper} maps from complete manifolds $Y$ to $X$, where ``{\it quasi-proper}'' means {\it locally constant at infinity.}
\end{description}

Clearly,
\begin{gather*}
 \operatorname{Sc}^{\rtimes_\downarrow}_\odot\geq \operatorname{Sc}^{\rtimes_\downarrow}\geq \operatorname{Sc}^{\rtimes_\downarrow} _{\widetilde {\rm sp}} \geq \operatorname{Sc}^{\rtimes_\downarrow} _{\rm sp}\geq
\operatorname{Sc}^{\rtimes_\downarrow} _{\rm st.par} \qquad \text{and}\qquad \operatorname{Sc}^{\rtimes_\downarrow} _{\widetilde{\rm sp}}\geq \operatorname{Sc}^{\rtimes_\downarrow} _{\pi_2=0}.
\end{gather*}

\textbf{II.} Assuming $X$ is a Riemannian manifold, relax the distance decreasing condition on $f$ by the following

$\bullet_{\rm area}$
the map $f$ {\it decreases the areas} of all surfaces in $Y$.\footnote{This makes sense for general metric spaces $X$
with the ``Hilbertian area'' defined in \cite{Gr2018'}.}

Clearly, $\operatorname{Sc}_{\rm area} ^{\rtimes_\downarrow}\geq \operatorname{Sc}^{\rtimes_\downarrow}$ and we show in \S2.E below that the ratio
 $\operatorname{Sc}_{\rm area} ^{\rtimes_\downarrow}\over \operatorname{Sc}^{\rtimes_\downarrow}$ can be arbitrarily large.

\textbf{III.} Replace the integer homology by the rational homology $H_\ast(X;\mathbb Q)$, which is essentially (but not quite)
the same as allowing maps $f\colon Y\to X$, where
$f_\ast$ sends
 the fundamental class of~$Y$ to a {\it non-zero multiple of $h$.}

\textbf{IV.} Instead of homology, use a bordism group of $X$, e.g., the spin bordism group, which is well adapted to $\operatorname{Sc}>0$.

\textbf{V.} {\it Remarks on singular $Y$.}

(a) If $h_m\in H_m(X)$ is not representable as $f_\ast[Y]$ for a smooth manifold $Y$, it may be interesting to try pseudomanifolds $Y$ with suitably defined singular Riemannian metrics~$g$  with ${\operatorname{Sc}(g)\geq \sigma}$.

(b) {\it Conical example.} Here $Y$ has an isolated singularity at a point $y_0\in Y$, where $g$ is a~smooth Riemannian metric on the complement to $y_0$ such that $\operatorname{Sc}(g)\geq \sigma$ and such that~$g $ is (approximately) {\it conical} at~$y_0$.

This means that there exists an $\varepsilon$-neighbourhood (ball) $U=U_\varepsilon\subset Y$ of~$y_0$, which topologically splits away from $y_0$,
\[
U\setminus \{y_0\}=Z\times (0,\varepsilon],
\]
 where $Z=(Z,g_Z)$ is a compact smooth Riemannian manifold such that the metric $g$ restricted to
 $U\setminus \{y_0\}$ is related to $g_Z$ as follows:
\[
g= a^2(t)t^2g_Z +{\rm d}t^2,
\]
where $a^2(t)>0$ is a smooth positive function on the (now closed) interval $[0, \varepsilon]$. (One may assume, if one wishes, that $a(t)$ is constant near $ t=0$.)

(c) One may additionally assume that $\operatorname{Sc}(g_Z)\geq \operatorname{Sc}\bigl(S^{m-1}\bigr)=(m-1)(m-2)$ for $m=\dim(Y)$ and, to make the metric truly conical, to require $a(t)$ to be constant near $ t=0$.

But this is not truly needed, since it can always be achieved by a small deformation of our $g$ near $y_0$.

(d) {\it Iterated conical singularities.} Next, following~\cite{Ch1979}, define $m$-dimensional (roughly) {\it cone-singular spaces} $Y$ with $\operatorname{Sc}(Y)\geq \sigma$ by induction on $m$, where (as in the above conical case) the metric (i.e., the distance function) on $Y$ is defined
by a smooth Riemannian metric $g$ on the non-singular part $Y_0\subset Y$, where the following conditions are satisfied:
\begin{enumerate}\itemsep=0pt
\item[(i)] The singular locus $\Sigma=Y\setminus Y_0$ is a closed subset in $Y$ with {\it codimension two} in $Y$.
\item[(ii)] $\operatorname{Sc}(g, y)\geq \sigma$ for all $y\in Y_0$.
\item[(iii)] Each $y_0\in Y$ admits a neighbourhood $U_0$, which is topologically (but not metrically) cylindrically splits away from $y_0$,
\[
U_0\setminus\{y_0\}=Z_0\times (0,\varepsilon_0],
\]
 where $Z_0=(Z_0, g_{Z_0})$, is a compact $(m-1)$-dimensional cone-singular space and where the Riemannian metric $g$ on the non-singular part
 of $U_0\setminus \{y_0\}$, denoted $U_{00} \subset U_0\setminus \{y_0\}$, is
\[
g= a_0^2(t)t^2g_{Z_0} +{\rm d}t^2+\delta_0,
\]
where $a_0^2(t)>0$ is a smooth positive function on $[0, \varepsilon_0]$ and where $\delta_0=\delta_0(u)$
 is a (small) smooth quadratic differential form on $U_0$, which {\it converges to zero} for $u\to y_0$.
\end{enumerate}

(e) Due to $h_0$, the above ``conical'' is slightly more general than how it is defined in (b) for an isolated singularity $y_0$.

(f) Similarly to the isolated singularity case, the requirement $\operatorname{Sc}(Z_0)\geq m-1)(m-2)$ does not significantly change the definition of $\operatorname{Sc}(Y)\geq \sigma$.

(g) One may also insist on the split-conical geometry at all points $y_0$:
{\it if $y_0$ is contained in the interior of an $l$-dimensional strata $S\subset \Sigma$, then
 a small neighbourhood $U_0\subset Y $ of $y_0$ metrically splits:
$U_0=S_0\times N_0$, where $S_0=U_0\cap S$, and where $N_0$ is a con-singular manifold with an $(m-l-1)$-dimensional base.}

(h) Probably, as in the isolated singularity case, this additional condition can be achieved by a small deformation of $g$ near $\Sigma$.

{\it Question.} How does the resulting
$\operatorname{Sc}^{\rtimes_\downarrow}(h)$, $h\in H_m(X)$, depend on the topology of the singular locus $\Sigma \subset Y$?

For instance, $Y$ may be iterated conical space with initial cones based on products of complex projective spaces and/or other generators of the oriented bordisms groups with $\operatorname{Sc}>0$ metrics (e.g., as in \cite{GL1983}), where we allow cones over $l$-dimensional $Y$ only if they admit metrics $g$ with ``suitably defined'' $\operatorname{Sc}(g)>0$ and/or, which is probably equivalent with $\operatorname{Sc}^\rtimes(g) >0$.

(Probably, stable minimal hypersurfaces and $\mu$-bubbles in such $Y$, similarly to how it is in a~smooth~$Y$, enjoy necessary properties required for the study of the scalar curvature
and this is also conceivable for the Dirac theoretic approach (compare with \cite{BR2021}).)

 \textbf {VI.} If $X$ is non-compact, allow classes $h$ with infinite supports\footnote{Sometimes referred to as ``locally finite homology classes'' as was pointed out to me by a referee.} and use proper (and quasi-proper) maps $f\colon Y\to X$.

{\it Remark on Completeness of $Y$.} Regardless of $X$ being complete or not, the value $\operatorname{Sc}_{\rm area} ^{\rtimes_\downarrow}[X]$ defined with complete $Y$ mapped to $X$ may be very different without this completeness.

For instance, $\operatorname{Sc}_{\rm area,sp} ^{\rtimes_\downarrow}[\mathbb R^n]$ defined with non-complete $Y$ is {\it infinite}.

 In fact, if $g_0$ is a metric on $\mathbb R^2$ such that $\operatorname{Sc}(g_0)=1$ and
 $\operatorname{area}(Y_0,g_0)=\infty$ (as in \S2.A below), then $Y=\bigl(\mathbb R^{n-2}\times \mathbb R^2, g_{\rm Eucl}+g_0\bigr)$ admits an obvious area contracting diffeomorphism onto $\bigl(\mathbb R^N,g_{\rm Eucl}\bigr)$.

But if we limit to {\it complete spin} manifolds $Y$, then
$\operatorname{Sc}_{\rm area,sp,compl} ^{\rtimes_\downarrow}[\mathbb R^n]=0$.

(It is unknown for $n\geq 4$ whether $\operatorname{Sc}_{\rm area,compl} ^{\rtimes_\downarrow}[\mathbb R^n]$ is zero or infinity without the spin assumption on~$Y$.)

\textbf{2.A. Surface examples.}
 Closed connected simply connected, i.e., {\it spherical} Riemann surfaces $X$ satisfy
\[
\operatorname{Sc}_{\rm area}^{\rtimes_\downarrow}[X]=\frac {8\pi} {\operatorname{area}(X)}.
\]

Indeed, the inequality $\operatorname{Sc}_{\rm area}^{\rtimes_\downarrow}[X]\geq\frac {8\pi} {\operatorname{area}(X)}$
follows from the existence of a measure preserving diffeomorphism from the 2-sphere with constant scalar curvature $\sigma=\frac{8\pi}{\operatorname{area}(X)}$ onto $X$;
 the opposite inequality follows from {\it Zhu's lemma} (see \cite{Zhu2019} and \cite[Section~2.8]{Gr2021}).

Similarly, one shows that
$\operatorname{Sc}_{\rm area}^{\rtimes_\downarrow}[X]$ for closed surfaces $X$ of positive genera.

On the opposite end of the spectrum, non-compact connected surfaces $X$ satisfy
$\operatorname{Sc}_{\rm area}^{\rtimes_\downarrow}[X]=\infty$,
since all the surfaces $X$ admit Riemannian metrics
with $\operatorname{Sc}=1$, and with the given areas (including ${\rm area}=\infty$ for non-compact~$X$) and since connected mutually diffeomorphic Riemann surfaces~$X_1$ and~$X_2$ of equal areas admit area preserving diffeomorphisms $X_1\leftrightarrow X_2$.

\textbf{Problem with $\boldsymbol{\operatorname{Sc}^{\rtimes_\downarrow}[X]=\operatorname{Sc}_{\rm dist}^{\rtimes_\downarrow}[X]}$.} Unlike $\operatorname{Sc}_{\rm area}^{\rtimes_\downarrow}$, the geometric meaning of
$\operatorname{Sc}^{\rtimes_\downarrow}[X]$ for spherical surfaces $X$ remains obscure. All one knows besides Zhu's lemma for general $X$ (see  \cite[Section~2.8]{Gr2021}) is that
\[
\operatorname{Sc}^{\rtimes_\downarrow}[X]<{4\pi^2\over \operatorname{diam}(X)^2}.
\]

\textbf{2.B. $\boldsymbol{H_{2m}(\mathbb CP^n)}$-example.} Let the complex projective space $\mathbb CP^n$ be endowed with
 the $U(n+1)$ invariant
(Fubini--Study)-metric such that the projective lines have scalar curvatures equal~2
 and let $\mathbb CP^m\hookrightarrow \mathbb CP^n$ be an $m$-plane.

If $m$ is odd, then the manifold $\mathbb CP^m$ is spin, and both homology and the spin-bordism class of
$ \mathbb CP^m \hookrightarrow \mathbb CP^n$ satisfy
\[
\operatorname{Sc}^{\rtimes_\downarrow}[\mathbb CP^m]\geq \operatorname{Sc}_{\rm sp}^{\rtimes_\downarrow}(k[\mathbb CP^m])\geq \operatorname{Sc}^{\rtimes_\downarrow}[\mathbb CP^m]_{\rm sp.brd}\geq \operatorname{Sc}( \mathbb CP^m)=m(m+1)
\]
and the same holds for the multiples $k[\mathbb CP^m]\in H_{2m}(\mathbb CP^n)$, $k=\dots,-1,0,1,2,\dots$.

If $k$ is even, then
if $m$ is {\it even}, then the $\operatorname{Sc}^{\rtimes_\downarrow}k[\mathbb CP^m]\geq m(m+1)$ remains valid for all $k$.
But if
 $k$ is {\it odd}, then
$\operatorname{Sc}_{\rm sp}^{\rtimes_\downarrow}(k[\mathbb CP^m])= 0$,
 since the classes
$k[\mathbb CP^m ]\in H_{2m}(\mathbb CP^n)$ are {\it not representable by the maps from spin manifolds} $Y\to \mathbb CP^n$.

In fact,
if an oriented manifold $Y^{2m} $ contains a smooth hypersurface~$H$ such that the $m$-fold self-intersection index $\underset {m}{\underbrace {H\frown \cdots \frown H}}$ is odd, then
the $(m-1)$-fold intersection is an orientable surface $\Sigma\subset Y$, which for even $m$ has
non-trivial normal bundle; hence $w_2[\Sigma]_{\mathbb Z_2}\neq 0$.

And if $k$ is {\it even}, then
\[
\operatorname{Sc}_{\rm sp}^{\rtimes_\downarrow}(k[\mathbb CP^m])\geq m^2
\]
since the class
$2[\mathbb CP^m])$ is represented by the
 quadric $Q^m\subset \mathbb CP^{m+1}\subset \mathbb CP^n$ given by the equation $z_0^2+z_1^2+\dots +z_m^2=0$,
 where this
$Q^m$ is spin and has scalar curvature~$m^2$.

Finally,
\[
\operatorname{Sc}_{\rm st.par}^{\rtimes_\downarrow}(m!h_{2m}) \geq {\rm const}_m \cdot \operatorname{Sc}\bigl(\underset {m} {\underbrace { S^2(1)\times \dots \times S^2(1)\bigr)}}=2m
\qquad \mbox{for all $m$ and $n\geq m$},
\]
since the quotient space
\[
\bigl(S^2\bigr)^m/\Pi(m) \qquad \text{of} \quad \bigl(S^2\bigr)^m=\underset {m} {\underbrace { S^2(1)\times \dots \times S^2(1)}}
\]
 by the permutation group $\Pi(m)$ admits a natural biholomorphic map $\psi\colon \bigl(S^2\bigr)^m\to \mathbb CP^m$,
 where ${\rm const}_m>0$ is the squared reciprocal to the minimal Lipschitz constant of
 maps in the homotopy class of this $\psi$.

 {\it Question.} What are $\operatorname{Sc}^{\rtimes_\downarrow}\bigl[\bigl(S^2\bigr)^m/\Pi(m)\bigr]$
 and of the symmetric powers $[(X)^m/\Pi(m)]$ for more general manifolds $X$?\footnote{These are among most attractive singular quasi-conical spaces discussed earlier.}

\textbf{2.C. Upper bounds and equalities.} The ($\mathbb T^\rtimes$-stabilized and $\widetilde{\rm sp}$-generalized) rigidity theorem by
Min-Oo \cite{Mi1995} and (the spin cobordims version of) Goette--Semmelmann's theorem from~\cite{GS1999} imply that the
class $h_{2m} =[\mathbb CP^m]\in H_{2m}(\mathbb CP^m)$ satisfy the following relations:
\begin{gather*}
\operatorname{Sc}_{{\rm area,sp.brd}\mathbb C P^m}^{\rtimes_\downarrow} (h_{2m})= m(m+1) \qquad \text{for all $m$ and $n\geq m$}, 
\end{gather*}
 where ``${\rm sp.brd}$'' indicates that this $\operatorname{Sc}_{\rm area}^{\rtimes_\downarrow}$ defined with smooth maps $Y\to X$ which are spin-bordant to the embedding $\mathbb CP^{m}\hookrightarrow \mathbb CP^{n}$,\footnote{It is unclear what happens for $m \leq n\leq 2m-2$.}
\begin{gather*}
\operatorname{Sc}_{{\rm area},\widetilde{\rm sp}}^{\rtimes_\downarrow} (h_{2m})= m(m+1) \qquad \text{for odd $m$},  
\\
\operatorname{Sc}^{\rtimes_\downarrow} _{{\rm area},\widetilde {\rm sp}}(h_{2m})=m^2 \qquad \text{for even $m$}, 
\\
\operatorname{Sc}^{\rtimes_\downarrow} _{{\rm area,st.par}}(h_{2m})_{\rm area} =2m \qquad \text{for all $m$ and $n\geq 2m-1$}. 
\end{gather*}

\textbf{2.D. Homological homogeneity conjecture.}
 Let $X$ be a compact symmetric space
 and $H\subset H_ m(X, \mathbb Q)$ be the linear subspace generated by the fundamental classes
 $[Y_i]\in H_m(X)$ of {\it homogeneous} (not necessarily totally geodesic) $m$-submanifolds $Y_i\subset X$.\footnote{$Y\subset X$ is {\it homogeneous} if an isometry group of $X$ preserves $Y$ and is transitive on $Y$.}
Then all classes
 ${h_m\in H}$
 can be represented by linear combinations of homogeneous~$Y_i$ such that
 $\operatorname{Sc}^{\rtimes_\downarrow} _{\rm area}[Y_i]\geq \operatorname{Sc}^{\rtimes_\downarrow}_{\rm area}(h)$.
(This maybe overoptimistic in general, but the $\widetilde{\rm sp}$-version of this can be, probably, proved with available means for products of spheres, complex and quaternionic projective spaces.)

\textbf{2.E. Equivalence conjecture.} All rational $h\in H_m(X) $ for compact Riemannian manifolds~$X$ without boundaries
satisfy:\footnote{The dimension $m=4$ may be special.}
\begin{gather*}
 \operatorname{Sc}^{\rtimes_\downarrow}_\odot(h)=
\operatorname{Sc}^{\rtimes_\downarrow}(h)\qquad \text{and} \qquad \operatorname{Sc}^{\rtimes_\downarrow}_{{\rm sp},\odot}(h)=
\operatorname{Sc}^{\rtimes_\downarrow}_{\rm sp}(h),\\
\operatorname{Sc}^{\rtimes_\downarrow}(h)\leq
 A\cdot \operatorname{Sc}^{\rtimes_\downarrow} _{\widetilde {\rm sp}}(h),\\
\operatorname{Sc}^{\rtimes_\downarrow} _{{{\widetilde sp } }}(h)\leq B\cdot \operatorname{Sc}^{\rtimes_\downarrow} _{\rm sp}(h),\\
\operatorname{Sc}^{\rtimes_\downarrow} _{\rm sp}(h)\leq
C\cdot \operatorname{Sc}^{\rtimes_\downarrow} _{\rm st.par}(h),
\end{gather*}
where $A=A_n$, $B=B_n$ and $C=C_n$ are universal constants,
and where the same relations are expected for the area version of these five $\operatorname{Sc}^{\rtimes_\downarrow}$.

 \textbf {2.F. Positivity.} Unlike $\operatorname{Sc}^\rtimes$, the values of all $\operatorname{Sc}^{\rtimes_\downarrow}=\operatorname{Sc}_{\rm dist}^{\rtimes_\downarrow}$-invariants
 are {\it non-negative}, since all (compact or not) manifolds $Y$ admit arbitrarily large Riemannian metrics with $\operatorname{Sc}>\varepsilon$.

 Moreover, the fundamental classes of {\it compact} connected manifolds $X$ with {\it non-empty boundaries}
 are {\it strictly positive}
 since such manifolds admit metrics with $\operatorname{Sc}>0$.
 (For instance, the $r$-balls with hyperbolic metrics $g$, ${\rm sect.curv}(g)=-1$, admit (obvious radial) metrics $g_+\geq g$ with $\operatorname{Sc}(g_+)\geq \exp-4r$.)

 \textbf {2.G. $\boldsymbol{[\exists \operatorname{Sc}>0]}$-Conjecture.} If a rational homology class $h\in H_m(X;\mathbb Q)$ vanishes under the
 classifying map $ \beta\colon X\to ${\sf B}$(\pi_1(X))$
\[
\beta_\ast (h)=0,
\]
 then $\operatorname{Sc}^{\rtimes_\downarrow}(h) >0$.\footnote{This seems more realistic if $\beta$ can be
homotoped to the $(m-2)$-skeleton of (some cell decomposition of) $ ${\sf B}$(\pi_1(X))$.}

\textbf {2.H. Finiteness.} $ \operatorname{Sc}^{\rtimes_\downarrow}(h)$ may be, a priori, infinite. However, the finiteness of $\operatorname{Sc}^{\rtimes_\downarrow}=\operatorname{Sc}^{\rtimes_\downarrow}_{\rm dist}$ easily follows from
 the $\square^m$-inequality~(3.8) in~\cite{Gr2021}, where the
 proof for $m\geq 9$ relies on Theorem~4.6 in~\cite{SY2017} and
 where the finiteness of $ \operatorname{Sc}_{\rm sp} ^{\rtimes_\downarrow}(h)$
 for all $m$ follows from~\cite{WXY2021}.
 (Probably, the arguments used in~\cite{WXY2021} generalize to $ \operatorname{Sc}_{\widetilde{\rm sp}} ^{\rtimes_\downarrow}$.)

 \textbf{2.I. $\boldsymbol{\nexists \operatorname{Sc}>0}$-Problem.}
 Does {\it non-vanishing} of a rational $h\in H_m(X;\mathbb Q)$ under the above classifying map {\sf B}$ (\pi_1(X))$
 imply that
$\operatorname{Sc} ^{\rtimes_\downarrow}=0$?
 This is known for $m=3$, and also for
 $\operatorname{Sc} ^{\rtimes_\downarrow}_{\widetilde{\rm sp}}(h)$ and all $m$ if the {\it spinorial curvature} $\mathbb S${\sf p}$.{\rm curv}^\downarrow(\beta_\ast(h)\in H_m(${\sf B}$(\pi_1(X))$
 defined in Section~\ref{section7} {\it vanishes}, e.g., if our {\sf B}$(\pi_1)$ admits a complete metric with ${\rm sect.curv}\leq 0$, see~\cite{Gr2021} and references therein.
 (I am not certain if there are examples of non-zero rational homology classes $\underline h$ in aspherical, say, compact finite-dimensional spaces such that $\mathbb S${\sf p}$.{\rm curv}^\downarrow(\underline h)\neq 0$.)

 \textbf {2.J. $\boldsymbol{\operatorname{Sc}_{\rm area}^{\rtimes_\downarrow}}$-finiteness question.} Is
$\operatorname{Sc}_{\rm area}^{\rtimes_\downarrow}[X]<\infty$
 for all compact Riemannian manifolds~$X$ without boundaries?

{\it Remarks.} (a) All metrics $g_+$ on a compact Riemannian
 manifold $(X,g)$ such that $\operatorname{area}_{g_+}(\Sigma) \geq \operatorname{area}_{g}(\Sigma)$ for all surfaces $\Sigma \subset X$ satisfy
\[
\operatorname{Sc}^\rtimes (g_+)\leq \mathrm{const}\cdot (X,g)<\infty.
\]
 In fact, this inequality holds for all $(Y, g_+)$ what admit area decreasing {\it spin maps}\footnote{A~continuous map between orientable manifolds, $f\colon Y\to X$, is \emph{spin} if $f^\ast(w_2(X))=w_2(Y)$, where $w_2$ is the second Stiefel--Whitney class.}
 $f\colon Y\to X$ with non-zero degrees.

 (b) Let $X=X_0\times Y$, where $Y$ is
 enlargeable\footnote{A compact Riemannian $n$-manifold $X$ is \emph{enlargeable} if there exists a sequence of oriented coverings
 $\tilde X_i\to X$ and {\it distance decreasing} maps
 $f_i\colon \tilde X_i\to S^n(R_i)$, $R_i\to \infty$, which are constant at infinity and which have {\it non-zero degrees}, compare with \cite{BH2009, GL1980, HKRS, SY1979}, in \cite[\S4.7]{Gr2021}, \cite[\S2.A]{Gr2023}.} and $\dim(X_0)=2$.
Then the finiteness of $\operatorname{Sc}_{\rm area}^{\rtimes_\downarrow}(X)$ for $X$ follows from \cite{Zhu2020}, where for $\dim(X)\geq 8$ one needs a version of Theorem~4.6 from \cite{SY2017}.

Also the
 $\rtimes$-stabilized version of the area slicing theorem from \cite{LM2021} (this stabilization is likely to be true) delivers an effective finite bound on $\operatorname{Sc}_{\rm area}^{\rtimes_\downarrow}(X_0\times Y)$.
For $\dim(X_0) =3$, provided $Y$ is enlargeable and $\dim(X)\leq 8$.

 But the principal case, where $X=S^n$ remains problematic for all $n\geq 4$
 and neither can one prove or disprove the existence of (necessarily non-spin) {\it complete} orientable $n$-manifolds $Y$, $n\geq 4$, with $\operatorname{Sc}\geq \sigma>0$,
which admit smooth proper area decreasing maps to $\mathbb R^n$, $n\geq 4$, with non-zero degrees.

\textbf{2.K. Outline of construction for $\boldsymbol{\operatorname{Sc}^{\rtimes_\downarrow}_{\rm area}/\operatorname{Sc}_{\rm dist}^{\rtimes_\downarrow}\to\infty}$.}
Let $g_0$ be a metric on a manifold~$Y$ such that $\operatorname{Sc}(g_0)>0$, then
{\it there exists metrics $g$ on $Y$ with arbitrarily large ratios $\operatorname{Sc}^{\rtimes\downarrow}_{\rm area}(g)/\operatorname{Sc}_{\rm dist}^{\rtimes\downarrow}(g)$.}
In fact, let $Y=(Y,g_0) $ be an arbitrary Riemannian manifold and $U\subset Y$ an open subset.
Then
for all $\varepsilon >0$ and $\delta>0$, there exists a Riemannian metric $g_{\varepsilon, \delta}$ on $Y$
 such that $ \operatorname{Sc}( g_{\varepsilon, \delta})\geq \operatorname{Sc}(g_0)-\varepsilon$ and
\begin{description}\itemsep=0pt\setlength{\leftskip}{0.50cm}
\item[$\bullet_{Y\setminus U}$] the metric $g_{\varepsilon, \delta}$ is equal to $g_0$ outside $U$;
\item[$\bullet_{\rm area}$] the metric $g_{\varepsilon, \delta}$ is area-wise smaller than $g_0$,
\[
\operatorname{area}_{g_{\varepsilon, \delta}}(S)\leq \operatorname{area}_{g_0}(S)\qquad \text{for all smooth surfaces $S\subset Y$};
\]
\item[$\bullet_{\rm dist}$] all Riemannian manifolds $X$, which 1-Lipschitz dominate\footnote{``Domination'' is a map with non-zero degree, see  \cite[\S1.5]{Gr2021}.}
 $(Y,g_{\varepsilon, \delta})$, have
$\operatorname{Sc}^\rtimes (X)\leq \delta$.
The only non-trivial condition here is $\bullet_{\rm dist}$, which is achieved with the follow ing
 one.
\item[$\bullet_{D}$] There is an open subset $U_D\subset U $ with $D=10\sqrt \frac {1}{\delta}$ such that
$(U_D, g_{\varepsilon, \delta})$ is isometric to the product
\[
\mathbb T^1(\varepsilon)\times \mathbb T^{n-2}(2\pi) \times [-D,D],
\]
where $T^1(\varepsilon)$ is the circle of length $\varepsilon$ and $\mathbb T^{n-2}(2\pi)$ is the standard flat torus.
\end{description}

The construction of $g_{\varepsilon, \delta}$, which satisfies $\bullet_{Y\setminus U}$,
$\bullet_{\rm area}$ and $\bullet_{D}$ is elementary (left to the reader\footnote{To get an insight, start with
$Y=S^2$, then look at $Y=S^2\times \mathbb T^{n-2}$.})
while the implication $\bullet_{D}\Rightarrow \bullet_{\rm dist}$ follows
the $\frac{2\pi}{n}$-inequality (see \cite[\S3.6]{Gr2021} and references therein).\footnote{The proof of $\frac{2\pi}{n}$-inequality for $n\geq 9$
relies on Theorem~4.6 in~\cite{SY2017}, and if one is satisfied with
$\operatorname{Sc}^{\rtimes\downarrow}_{\rm area}(g)/\operatorname{Sc}_{{\rm dist,sp}}^{\rtimes\downarrow}(g)\to \infty$, then one can use the spinorial version of
$\frac{2\pi}{n}$ from~\cite{Ze2019}.}

Thus, we see that
 {\it if a non-torsion homology class $h$ in a compact manifold $X$ satisfies
 $\operatorname{Sc}_{\rm dist}^{\rtimes_\downarrow}(h)>0$, then the ratio
 $\operatorname{Sc}^{\rtimes_\downarrow}_{{\rm area,sp}}(h)/\operatorname{Sc}_{{\rm dist,sp}}(h)^{\rtimes_\downarrow}$
 can be made arbitrarily large with some Riemannian metric on~$X$.}

\textbf{2.L. Question on $\boldsymbol{\lambda^\downarrow_1(h,\beta)}$.} The definition
of $\operatorname{Sc}^{\rtimes_\downarrow}(h)$, which depends on $\lambda_1(Y,\beta=1/4)$
(see $(N + 1)/4N$-Remark 1.D) makes sense for all $\beta$ and the arguments which depends on stable minimal hypersurfaces and $\mu$-bubbles generalize to all $\beta$, e.g., as the $\square^{\exists \exists}(n, m,N)$-inequality, which is stated in \cite[\S2.B]{Gr2023} for $\beta=N/(N+1)$. However, the geometric significance of this for
$\beta\neq N/(N+1)$ is unclear.

Probably, if $X$ is simply connected, $\beta \leq \beta_m>0$ and $m\geq 3$, then an integer multiples
$lh$ for some $l\neq 0$ and all $h\in H_m(X)$ are representable by a~distance decreasing maps $Y\to X$, where $\lambda_1(Y, \beta)\geq C$ for a given $C>0$.

\subsection*{Exercises}

 \textbf{2.M.} Let $g$
be a Riemannian metric on an {\it open manifold}\footnote{A manifold $X$ is \emph{open} if it contains no closed manifold connected component.} $X$ of dimension $\dim(X)=n\geq 2$. Show that there exists a Riemannian metric $g_+$ on $X$ such that $\operatorname{Sc}(g_+)=1$ and $\operatorname{area}_{g_+} (\Sigma)\geq \operatorname{area}_{g} (\Sigma)$ for all smooth surfaces $\Sigma\subset X$.

{\it Hint}: Observe that $[0,1]\times \mathbb R^{n-1}$ admits an area decreasing diffeomorphism onto
$\mathbb R^n$ and
use products of surfaces with constant curvatures by $\mathbb R^{n-2}$ as building blocks for $(X,g_+)$.

 {\it Remark.} If $Y$ is a {\it complete} spin $n$-manifold with $\operatorname{Sc}^\rtimes (Y) \geq \sigma>0$, then it admits no proper area decreasing map to $\mathbb R^n $ with non-zero degree~\cite{GL1983}.

 \textbf {2.N.} Show that non-zero multiples of homology classes $h$ in {\it simply connected} manifolds $X$
 have $\operatorname{Sc}_{{\rm st.par}}^{\rtimes_\downarrow}(ih)>0$, for some $i\neq 0$.

 {\it Hint.} Recall the Serre--Thom theorem on framed bordisms
 and apply Stolz' theorem on spin manifolds~\cite{St1992}.

\section[rtimes\^{}\{downarrow\}-extremality and rtimes\^{}\{downarrow\}-rigidity]{$\boldsymbol{\rtimes^\downarrow}$-extremality and $\boldsymbol{\rtimes^\downarrow}$-rigidity}\label{section3}

\textbf{3.A. Homological $\boldsymbol{\operatorname{Sc}_\ast^{\rtimes_\downarrow}}$-problems.} Let $X$ be a Riemannian manifold and $h\in H_m(X)$ a~homology class,
e.g., $m=n=\dim(X)$, and let $h$ be the fundamental class $[X]$ of $X$, where $X$ is assumed oriented.

Evaluate $\operatorname{Sc}_\ast^{\rtimes_\downarrow}$ and/or
find relations between $\operatorname{Sc}_\ast^{\rtimes_\downarrow}$
 and more accessible metric invariants of~$X$.

Decide if $\operatorname{Sc}_\ast^{\rtimes_\downarrow}(h)$ is represented
by an {\it $\operatorname{Sc}_\ast^{\rtimes_\downarrow}$-extremal}, or, for brevity, {\it $\rtimes^\downarrow_\ast$-extremal},
oriented Riemannian $m$-manifold mapped to $X$,
\[
Y\overset {f} \rightarrow X \quad \text{such that $f_\ast(Y)=h$ and $\operatorname{Sc}^\rtimes (Y)= \operatorname{Sc}_\ast^{\rtimes_\downarrow}(h)$},
\]
where $f$ is the distance or the area decreasing depending on ``$\ast$'' and where, ideally, $f$ is an isometric immersion.

For instance, given a submanifold $Y\hookrightarrow X$, e.g., $Y=X$ decide if it is $\operatorname{Sc}_\ast^{\rtimes_\downarrow}$-extremal,
or, moreover, if it is {\it rigid}, that is {\it unique extremal} (compare with \S3.D below).

Find examples of $h$, where there is no extremal manifold $Y\to X$ with $f_\ast[Y]=h$, but such a generalized $Y$, e.g., a singular extremal one does exist. (We saw some potential examples of such singular $Y$, and stable minimal singular hypersurfaces suggest further examples.)

Determine which closed manifolds $X$ admit
$\rtimes^\downarrow$-extremal Riemannian metrics.

(Possibly, metrics $g_0$ with $\operatorname{Ricci}(g_0) >0$ can be deformed to $g\geq g_0$ with
$\operatorname{Sc}^{\rtimes}(g)= \operatorname{Sc}^{\rtimes_\downarrow}(g)$.
But, for instance, metrics $g_0=g_1+g_2$ on $X = X_1\times X_2 $,
 where ${\rm sect.curv}(g_1)=1$ and ${\rm sect.curv}(g_2)< 0$ admit no such deformations. However, the pointed Hausdorff limit manifolds $\lim_{\lambda\to \infty}(X,g_+\lambda g_2)$, which are isometric to $X_1\times \mathbb R^{\dim(X)}$,
 are extremal.

In general, the existence of a metric $g_0$ on $X$ with $\operatorname{Ricci}(g_0)\geq 0$ might be necessary for the existence of an extremal metric $g$ on $X$.

\subsection*{Examples}

\textbf{3.B$\boldsymbol{{}_{S^n}}$.}
 Complete manifolds with {\it constant sectional curvatures}, e.g., unit spheres, flat tori and Euclidean spaces
 are {\it $\operatorname{Sc}^{\rtimes_\downarrow}_{{\rm area},\widetilde{\rm sp}}$-extremal}.

 This follows from the $\rtimes$-stabilized {\it Llarull's theorem} (see~\cite{Gr2021} and references therein).

\textbf{3.B$\boldsymbol{{}_{\mathcal R>0}}$.} A compact spin manifold $X$
 with {\it non-negative} curvature operator, $ \mathcal R(X)\geq 0$,
e.g., a compact symmetric space is $\operatorname{Sc}^{\rtimes_\downarrow}_{{\rm area},\widetilde{\rm sp}}$-extremal, provided {\it scalar curvature $\operatorname{Sc}(X)$ is constant}\footnote{There are lots of metrics with $\mathcal R>0$ on spheres $S^n$ and if $n\geq 3$ many of these have constant scalar curvatures.
On the other
 hand, it is possible that a closer look at the curvature term in the twisted Schr\"odinger--Lichnerowicz formula (see Section~\ref{section5}) would allow one to drop the constancy of the scalar curvature condition.}
 and {\it the Euler characteristic of the universal covering $\tilde X$ does not vanish}.

 This follows
 by an elaboration on the proof of the {\it Goette--Semmelmann extremality theorem}~\cite{GS2002}. (We say a few words about it in Section~\ref{section5}.)

Probably, the corresponding rigidity arguments (see~\cite{Li2010} and references therein) also
 admit $\rtimes $-stabilization, but I did not check this carefully.

Also the
 condition $\chi\bigl(\tilde X\bigr)\neq 0$ seems redundant and $\rtimes^{\downarrow}_{{\rm area},\widetilde{\rm sp}}$-extremality can be strengthened, also conjecturally, to
 the $\rtimes^{\downarrow}_{\rm area}$-extremality.

\textbf{3.B$\boldsymbol{{}_{\bigtimes [a_1,b_i]}}$.} The rectangular solids $\bigtimes_1^n[-a_i,b_i]$
 are $\rtimes^\downarrow_{\rm sp}$-extremal and, if $n\leq 8$, they are ${\rtimes^\downarrow}$-extremal.

 In fact, $\rtimes^{\downarrow}_{\rm sp}$-extremality follows by a slight generalization of the argument from \cite{WXY2021}, which,
probably, can be adapted for the proof of the
 $\rtimes^{\downarrow}_{\widetilde{\rm sp}}$-extremality.

 As for $\rtimes^{\downarrow}$-extremality for $n\leq 8$, this follows from
 the $\square^{\exists\exists}(n,m,N)$-inequality  \cite[\S2.B]{Gr2023}.

 Furthermore, the generic regularity theorem from \cite{CMS2023} extended to $\mu$-bubbles (I have not check this extension) yields the $\square^{\exists\exists}(n,m,N)$-inequality and thus $\rtimes^{\downarrow}$-extremality of solids for $n\leq10$.

 Moreover, granted
 a $\mu$-bubble generalization of Theorem 4.6 from~\cite{SY2017}, the $\rtimes^{\downarrow}$-extremality (but not the $\square^{\exists\exists}(n,m,N)$-inequality) would follow for all~$n$.

 \textbf{3.B$\boldsymbol{{}_{\bigtimes \bigtimes}}$.}
 Riemannian products of the manifolds from the above examples,
 e.g.,
 \[
 X= \bigl(\bigtimes_1^{n-k}[-a_i,b_i]\bigr)\times S^k,
 \]
  are $\rtimes^{\downarrow}_{\rm sp}$-extremal.

As above, this follows by a simple generalization of argument from \cite{WXY2021}
 combined with the basic (algebraic) inequality in \cite{GS2002} for twisted Dirac operators on manifolds
 with $\mathcal R\geq 0$.

But the $\rtimes^{\downarrow}_{\rm sp}$-extremality remains problematic even for $n\leq 8$.

 For instance, if $k\leq 4$ and $n\leq 8$ (probably $n\leq 10$ will do), then the $\square^{\exists\exists}(n,m,N)$-inequality combined with the warped product splitting argument in \cite[\S5.5]{Gr2021} yield
 $\rtimes^{\downarrow}_{\rm sp}$-extremality of $X= \bigl(\bigtimes_1^{n-k}[-a_i,b_i]\bigr)\times S^k$.

 Yet, there is no approach so far to non-spin extremality of the spheres $S^k$ for $k\geq 5$.\footnote{The warped product splitting argument (combined with a stable version of~\cite{GS2002}) applies to $S^4$, because 3-manifolds are spin.}

 \textbf{3.B$\boldsymbol{{}_{\rm warp}}$.} There are several classes of log-concave
 warped product manifolds,
 e.g., $S^n$ minus a~point, where the $\rtimes_{\rm sp}$-extremality (and $\rtimes$-extremality for $n=4$) follow by
 $\rtimes$-stabilization of the arguments in \cite[\S\S5.5--5.7]{Gr2021} and \cite{CZ2021a}.
  In fact, the $\rtimes$-extremality for warped manifolds is more common then non-stabilized extremality.

 For instance,
 geodesic balls in spheres and in $\mathbb R^n$ are not non-stably extremal:
 one can increase their metrics without diminishing the scalar curvatures. But,
probably, they are
 $\rtimes_\downarrow$-extremal.\looseness=1

 \textbf{3.C. Questions.}
 \begin{enumerate}\itemsep=0pt
 \item[(i)] Which convex subsets in $\mathbb R^n$ are $\rtimes^\downarrow$-extremal?
\item[(ii)] Which surfaces are $\rtimes^\downarrow$-extremal?
\end{enumerate}

\textbf{3.D. About rigidity.}
 The proofs of extremality of the manifolds $X$ in the above examples can be upgraded to {\it rigidity} that says in the present case that
 {\it if a smooth distance non-increasing positive degree map $f\colon Y\to X$ satisfies
 $\operatorname{Sc}_\ast^\rtimes(Y)\geq \operatorname{Sc}_\ast^{\rtimes\downarrow}(X)$ $($where $ \operatorname{Sc}_\ast^{\rtimes\downarrow}(X)=\operatorname{Sc}_\ast^\rtimes(X)$ by extremality$)$, then $f$ is homotopic to a local isometry, where one
 can drop ``homotopic to''
 if $X$ has no local scalar flat factors.}

This follows by combining the $\rtimes$-stabilized rigidity arguments in \cite{GS2002} and \cite{Li2010} with those in \cite[\S5.7]{Gr2021}
but to be honest, I did not check this in full generality.

 \section[Sc\^{}\{rtimes downarrow\}-product inequalities, conjectures and problems]{$\boldsymbol{\operatorname{Sc}^{\rtimes\downarrow}}$-product inequalities, conjectures and problems}\label{section4}

 \textbf{4.A. Additivity for cylinders.} Since, obviously,
\[
\operatorname{Sc}^{\rtimes_\downarrow}[X_1\times X_2]\geq \operatorname{Sc}^{\rtimes_\downarrow}[X_1]+\operatorname{Sc}^{\rtimes_\downarrow}[X_2],
\]
then, for all Riemannian manifolds $X_1$ and $X_2$, the inequality
\[
\operatorname{Sc}^{\rtimes_\downarrow}[X_1\times X_2]\leq \operatorname{Sc}^{\rtimes_\downarrow}[X_1]+\operatorname{Sc}^{\rtimes_\downarrow}[X_2],
\]
 is equivalent to the equality
\[
\operatorname{Sc}^{\rtimes_\downarrow}[X_1\times X_2]= \operatorname{Sc}^{\rtimes_\downarrow}[X_1]+\operatorname{Sc}^{\rtimes_\downarrow}[X_2].
\]

Thus, in particular,
 the $\square^{\exists\exists}(n,m,N)$-{\it inequality} from \cite{Gr2023} and/or {\it equivariant separation theorem}
for stable $\mu$-bubbles\footnote{See \cite[\S5.4]{Gr2021} and compare with~\cite{GZ2021, Ri2020} and with
 \cite[the proof of \S2.B]{Gr2023}.} along with the equality
\[
\operatorname{Sc}^\rtimes[a,b]=4\lambda_1[a,b]= \frac {4\pi^2}{(b-a)^2}.
\]
 imply the following.

{\it Proposition.} The fundamental homology classes of oriented Riemannian cylindrical manifolds $X=Y\times [a,b]$ of dimensions $\leq 8$ satisfy
\[
\operatorname{Sc}^{\rtimes_\downarrow}[X]= \operatorname{Sc}^{\rtimes_\downarrow}[Y]+\operatorname{Sc}^{\rtimes_\downarrow}[a,b].
\]
 (This generalizes $\operatorname{Sc}^{\rtimes_\downarrow}\left( \bigtimes_1^n[a_i,b_i]\right)=
\sum_i\operatorname{Sc}^{\rtimes\downarrow}[a_i,b_i]$, that is \S3.B$_{\bigtimes [a_1,b_i]}$ from the previous section.)

\textbf{4.B. The spin case.} This additivity formula remains {\it problematic} for
 $n\geq 9$,\footnote{The dimensions $n=9,10$, probably, can be taken care by the argument in~\cite{CMS2023}.} but {\it the spin cube inequality} from \cite{WXY2021} (proved with an index theorem for deformed Dirac operators on manifolds with boundaries) implies, as we stated earlier, that{\samepage
\[
\operatorname{Sc}_{\rm sp}^{\rtimes\downarrow} ( \bigtimes_1^n[-a_i,b_i] )=\operatorname{Sc}^{\rtimes} ( \bigtimes_1^n[-a_i,b_i] )=
\sum_1^n \frac {4\pi^2}{(b_i-a_i)^2} =
\sum_1^n \operatorname{Sc}^{\rtimes_\downarrow}_{\rm sp}[a_i,b_i]
\]
for all $n$.}

Yet, as far as I can see, the present day Dirac theoretic argument does not yield
the general $\operatorname{Sc}_{\rm sp}$-inequality
\[
\operatorname{Sc}_{\rm sp}^{\rtimes_\downarrow}[Y\times [a,b] ]\leq \operatorname{Sc}_{\rm sp}^{\rtimes_\downarrow}[Y]+\operatorname{Sc}_{\rm sp}
^{\rtimes_\downarrow}[a,b].
\]
 However, this argument does apply, if $Y$ is a special (extremal) manifold as in \S3.B$_{\times \times}$, e.g., a~product of spheres.

 \textbf {4.C. [$\bf sect.curv\leq 0$]-Remark.} Let $Y$ and $Z$ be compact Riemannian manifolds, where $Z$ has no boundary and the sectional
 curvature ${\rm sect.curv}(Y)\leq 0$. Then, similarly as above, one can prove additivity in the following two cases:
\begin{enumerate}\itemsep=0pt
\item[(i)] If $\dim(Y\times Z)=n\leq 8$,\footnote{In view of \cite{CMS2023}, the inequality $n\leq 10 $ may suffice.} then
\[
\operatorname{Sc}^{\rtimes_\downarrow}[Y\times Z]=\operatorname{Sc}^{\rtimes_\downarrow}[Y]=\operatorname{Sc}^{\rtimes_\downarrow}[Y] +\operatorname{Sc}^{\rtimes_\downarrow}[Z]=\operatorname{Sc}^{\rtimes_\downarrow}[Y]. \]
\item[(ii)] If $Y$ is as in \S3.B$_{\times \times}$, then
 \[
 \operatorname{Sc}_{\rm sp}^{\rtimes_\downarrow}[Y\times Z]= \operatorname{Sc}_{\rm sp}^{\rtimes_\downarrow}[Y]=\operatorname{Sc}_{\rm sp}^{\rtimes_\downarrow}[Y]+\operatorname{Sc}_{\rm sp}^{\rtimes_\downarrow}[Z]=\operatorname{Sc}^{\rtimes_\downarrow}[Y]
 \]
for all $n$.
\end{enumerate}

\textbf {4.D. Riemannian additivity conjecture.}
Riemannian products of all oriented Riemannian manifolds satisfy
\[
\operatorname{Sc}^{\rtimes_\downarrow}[X_1\times X_2]= \operatorname{Sc}^{\rtimes_\downarrow}[X_1]+ \operatorname{Sc}^{\rtimes_\downarrow}[X_2].
\]

In fact, the following stronger inequality might be true.

\textbf {4.E.} \textbf{Sup-metric product conjecture.} Let $X_i$, $i=1,\dots,k$, be metric spaces (e.g., closed oriented Riemannian
 manifolds) and let
\[
X=\bigtimes_{\sup_i} X_i=
 (X_1
 \times \cdots \times X_k, \operatorname{dist}_{\sup})
 \]
be their product endowed with the {\it sup-metric}
\[
\operatorname{dist}((x_1,\dots,x_k), (x'_1,\dots,x'_k))=\max_{i=1,\dots,k} \operatorname{dist}(x_i, x'_i).
\]
Then rational homology classes $h_i\in H_{m_i}(X_i;\mathbb Q)$ (e.g., the rational fundamental classes $[X_i]$\footnote{``Rational'' in the case of compact locally contractible spaces means ``{\it a non-zero integer multiple of}'', that is,
\[
\operatorname{Sc}^{\rtimes_\downarrow}( h_\mathbb Q) \overset{\rm def}{=}
\sup_{N\neq 0}\operatorname{Sc}^{\rtimes_\downarrow}(N\cdot h_\mathbb Q).
\]})
satisfy
\begin{equation}\label{eq-b}
 \operatorname{Sc}^{\rtimes_\downarrow}( \otimes_i h_i)\leq\sum_{i=1,\dots,k} \operatorname{Sc}^{\rtimes_\downarrow}( h_i ), \qquad \text{e.g.,} \qquad \operatorname{Sc}^{\rtimes_\downarrow}[ \bigtimes_{\sup_i} X_i]_{\mathbb Q}\leq\sum_{i=1,\dots,k} \operatorname{Sc}^{\rtimes_\downarrow} [X_i]_{\mathbb Q},
 \end{equation}
where the opposite inequality
\[
\operatorname{Sc}^{\rtimes_\downarrow}( \otimes_i h_i )\geq\sum_{i=1,\dots,k} \operatorname{Sc}^{\rtimes_\downarrow}( h_i )
\]
follows from additivity of the scalar curvature;
hence, \eqref{eq-b} implies the equality
\[
 \operatorname{Sc}^{\rtimes_\downarrow} ( \otimes_i h_i) =\sum_{i=1,\dots,k} \operatorname{Sc}^{\rtimes_\downarrow} (h_i).
 \]

\textbf{4.F. $\boldsymbol{\bigtimes_{\sup_i}[a_i,b_i]}$-Example.}
The above indicated proofs of \S4.A and \S4.B actually show that the rectangular solids $\bigtimes_1^n[a_i,b_i]$ with the Riemannian product and the sup-product metrics have the same $\operatorname{Sc}^{\rtimes_\downarrow}$ for $n\leq8$ and
have the same $\operatorname{Sc}^{\rtimes_\downarrow}_{\rm sp}$ for all~$n$.
This {\it confirms the validity of~\eqref{eq-b} for rectangular solids}.

\textbf{4.F$'$. $\boldsymbol{\bigtimes [0,d_i]}$-Sub-Example.} Let $Y\subset\mathbb R^n$ be a diffeomorphic image
 of the $n$-cube and let~$d_i$, $ i=1,\dots,n$, be the distances between
 the images in $Y$ of the pairs of the opposite $(n-1)$-faces of the cube.
Then the first Dirichlet eigenvalue of the Laplacian $-\Delta_Y$
is bounded by that of the solid $\bigtimes_i [0,d_i]$,
\[
\lambda_1(-\Delta_Y)\leq \sum_i\frac{\pi^2}{d^2_i}.
\]

{\it Exercise.} Find a direct elementary proof of this inequality.\footnote{To my shame, I could not solve it.}

\textbf{Sup-distance, sup-area and $\boldsymbol{\operatorname{Sc}_{\rm sup.area}^{\rtimes_\downarrow}}$.} The Riemannian product metric, that is the Pythagorean one
\[
\operatorname{dist}((x_1,\dots,x_k), (x'_1,\dots,x'_k))=\sqrt{ \sum_i\operatorname{dist}(x_i, x'_i)^2},
\]
 is greater than the sup-metric but only by a factor $ \sqrt k$,
\[
1\leq {\operatorname{dist}((x_1,\dots,x_k), (x'_1,\dots,x'_k))\over \operatorname{dist}_{\sup}((x_1,\dots,x_k), (x'_1,\dots,x'_k))}\leq \sqrt k.
\]

The situation is somewhat different with areas. Namely, let $X=\bigtimes_iX_i$
 be the product of Riemannian manifolds and let
 $\sup_i$-$\operatorname{area}(\Sigma)$ for a smooth surface $\Sigma\subset X$ be the maximum of the areas of the projections $\Sigma\to X_i$.
 Here again
 \[
 {\sup_i}\mbox{-}{\operatorname{area}}(\Sigma)\leq \operatorname{area}(\Sigma)
 \]
 but now, unlike to how it is with the distances, the ratio
\[
\frac{\operatorname{area}(\Sigma)}{{\sup_i}\mbox{-}{\operatorname{area}}(\Sigma)}
\]
 may be infinite.
 Accordingly, the corresponding $\operatorname{Sc}_{{\rm sup.area}}^{\rtimes_\downarrow}(h)$,
 $h\in H_m(
 \bigtimes_iX_i)$, defined with smooth maps $f\colon Y^m\to \bigtimes_iX_i$,
 $f_\ast[Y]=h$ such that the corresponding $f_i\colon Y^m\to X_i$ are area decreasing, can be {\it significantly greater}
 than $\operatorname{Sc}_{{\rm sup.area}}^{\rtimes_\downarrow}(h)$, where the maps $f$ must be area decreasing themselves.

 Thus, the area version of \eqref{eq-b},
\begin{equation}\label{eq-c}
\operatorname{Sc}_{\rm sup.area}^{\rtimes_\downarrow}( \otimes_i h_i)\leq\sum_{i=1,\dots,k} \operatorname{Sc}_{\rm area}^{\rtimes_\downarrow}( h_i ),
\end{equation}
 e.g.,
 \[
 \operatorname{Sc}_{\rm sup.area}^{\rtimes_\downarrow}[ \bigtimes_{\sup_i} X_i]_{\mathbb Q}\leq\sum_{i=1,\dots,k} \operatorname{Sc}_{\rm area}^{\rtimes_\downarrow} [X_i]_{\mathbb Q},
 \]
is qualitatively stronger than corresponding inequality for $\operatorname{Sc}_{\rm area}^{\rtimes_\downarrow}( \otimes_i h_i).$

Although we have no known means for bounding $\operatorname{Sc}_{\rm area}^{\rtimes_\downarrow}$ and even less for $\operatorname{Sc}_{{\rm sup.area}}^{\rtimes_\downarrow}$ in most cases, we shall
do this in the next section for $\operatorname{Sc}_{{\rm sup.area}, \widetilde{\rm sp}}^{\rtimes_\downarrow}$ and thus prove
the $\widetilde{\rm sp}$-version of~\eqref{eq-c}
in some cases.

\textbf{4.G. Semiadditivity problem.} Let $ X=X^n$ and $Z=Z^k$, $k\leq n-2$, be compact Riemannian manifolds, possibly with boundaries,
and let $f\colon X\to Z$ be a smooth distance decreasing map such that $\partial X\overset {f}\to\partial Z$,
and let $h_m=\bigl[f^{-1}(z)\bigr]\in H_m(X)$, $m=n-k$, be the homology class of the pullback of a generic $z\in Z$.

Identify the cases, where
\[
\operatorname{Sc}^{\rtimes\downarrow}(h_m)_{\mathbb Q}\geq \operatorname{Sc}^{\rtimes\downarrow}[X]_{\mathbb Q} - \operatorname{Sc}^{\rtimes\downarrow}[Y]_{\mathbb Q} ,
\]
at least for ``simple'' manifolds $Z$, e.g.,
compact convex domains in $\mathbb R^k$ and in $S^k$ and, in general,
evaluate the difference
\[
\operatorname{Sc}^{\rtimes\downarrow}[X]_{\mathbb Q} -\operatorname{Sc}^{\rtimes\downarrow}[Y]_{\mathbb Q} -\operatorname{Sc}^{\rtimes\downarrow}(h_m)_{\mathbb Q}
\]
in terms of the geometry of $Z$, for instance, where $Z$ is the product of balls $Z=\bigtimes_iB^{k_i}(R_i)$ or product of spheres $S^{k_i}(R_i)$.

If $n=m+k\leq 8$, a satisfactory lower bound on $\operatorname{Sc}^{\rtimes\downarrow}(h_m)$
 for rectangular solids $Z$ follows from~\S4.A.
 Also \cite[\S2.B]{Gr2023} yields similar bounds for products of 2-discs and 2-spheres
(compare~\cite{GZ2021}). But it is unclear, for instance, how large the difference $\operatorname{Sc}^{\rtimes\downarrow}[X]_{\mathbb Q} -\operatorname{Sc}^{\rtimes\downarrow}[Y]_{\mathbb Q} -\operatorname{Sc}^{\rtimes\downarrow}(h_m)_{\mathbb Q}$
can be for the balls $B^k\subset \mathbb R^k$ and spheres $S^k$ for large~$k$.

\section{Additivity of the twisted SLWB-formula and applications}\label{section5}

Let $Y$ be a Riemannian spin $n$-manifold and $V\to X$ be a complex vector bundle with a unitary connection $\nabla$ and let $\mathcal D_{\otimes V}$ denote
the Dirac operator on spinors $\mathbb S$ on $Y$ tensored with $V$. Then the square of $\mathcal D_{\otimes V}$ satisfies the following
 {\it Schr\"odinger--Lichnerowicz--Weitzenb\"ock--Bochner formula} (see~\cite{LM1989})
\[
{\cal D}_{\otimes V}^2=\nabla^2_{\otimes V} + \frac{1}{4}\operatorname{Sc}(Y) + {\cal K}_{\otimes V},
\]
where $\mathcal K_{\otimes V}$ is an endomorphism
$\mathbb S\otimes V\to \mathbb S\otimes V $ such that
\[
\mathcal K_{\otimes V}(s\otimes v)=\frac {1}{2}\sum_{i,j}(e_i\circ e_j\circ s) \otimes R^V_{e_i\wedge e_j}(v),
\]
where $e_i\in T_y(Y) $, $i=1,\dots,m=\dim(Y)$, are orthonormal tangent vectors at $y\in Y$, where
 $\circ$ is the Clifford multiplication and $R^V_{e_i\wedge e_j}\colon V\to V$ is the curvature operator of $\nabla$.

Next, recall that the curvature of the tensor product of two bundles with connections satisfies
\[
R^{V_1\otimes V_2}=1^{V_1}\otimes R^{ V_2}+R^{V_1}\otimes 1^{V_2},
\]
where $1^V\colon L\to V$ is the identity operator, and observe that
 the operators on $\mathbb S\otimes V_1\otimes V_2$ defined by
\[
s\otimes v_1\otimes v_2\mapsto \frac {1}{2}\sum_{i,j}(e_i\circ e_j\circ s)
 \otimes R^{V_1}_{e_i\wedge e_j}(v_1)\otimes v_2
 \]
and by
\[
s\otimes v_1\otimes v_2\mapsto \frac {1}{2}\sum_{i,j}(e_i\circ e_j\circ s)
 \otimes v_1\otimes R^{V_2}_{e_i\wedge e_j}(v_2)
 \]
have the same spectra up to multiplicity as
\[
s\otimes v_1\mapsto \frac {1}{2}\sum_{i,j}(e_i\circ e_j\circ s) \otimes R^{V_1}_{e_i\wedge e_j}(v_1)
\]
and
\[
s\otimes v_2\mapsto \frac {1}{2}\sum_{i,j}(e_i\circ e_j\circ s)
 \otimes R^{V_2}_{e_i\wedge e_j}(v_2)
\]
correspondingly.
Therefore, the lowest eigenvalue $\lambda_{\otimes1\otimes 2} $ (often negative) of the (self-adjoint) operator $\mathcal R_{\otimes (V_1\otimes V_2)}$ is bounded from below by the sum of these for
$\mathcal R_{\otimes V_1}$ and $\mathcal R_{\otimes V_2}$,\footnote{A referee pointed out to me that $\lambda_{\otimes 1\otimes 2}= \lambda_{\otimes 1}+\lambda_{\otimes 2}$.}
\[
\lambda_{\otimes 1\otimes 2}\geq \lambda_{\otimes 1}+\lambda_{\otimes 2}.
\]

This yields the following.

\textbf {5.A. Theorem.}\footnote{This is a refinement of the Llarull--Goette--Semmelmann--Listing rigidity theorem.}
Let $X=\bigtimes_k X_k$, $k=1,\dots,l$, be an orientable Riemannian $n$-manifold split into Riemannian product, where the factors $X_k =(X_k,\underline g_k)$ are
either
\begin{enumerate}\itemsep=0pt
\item[(a)] compact $n_k$-manifolds with {\it non-negative curvature operators}, $\mathcal R^{X_k}\geq 0$ (e.g., closed convex hypersurfaces in $\mathbb R^{n_k+1}$)
 and with {\it non-vanishing Euler characteristics} $\chi(X_k)\neq 0$ (hence of even dimensions $n_k$),
 or
\item[(b)] spheres $S^{n_k}$ with constant sectional curvatures $($possibly of odd dimension $n_k)$.
\end{enumerate}

 Let $\underline g^\natural_k =\operatorname{Sc}(\underline g_k)\cdot \underline g_k$.
 Let $Y=(Y,g) $ be a smooth complete orientable Riemannian $(n+N)$-manifold with $\operatorname{Sc}(g) >0$,
and let $g^\natural =\operatorname{Sc}(g) \cdot g$.
Let
 $Z$ be an orientable enlargeable $N$-manifold, e.g., $Z=\mathbb R^N$,
and let
$f\colon Y\to X\times Z$ be a~smooth proper (quasi-proper will do) map such that
the corresponding maps $f_k\colon Y\to X_k$ are {\it strictly sum-wise area decreasing} with respect to $g^\natural$
in~$Y$ and~$\underline g^\natural_k $ in~$X_k$.\footnote{The quadratic forms $\underline g^\natural_k$ on $ X_k$ may vanish but ``area'' makes sense anyway.}

This means that the norms of the exterior squares of the differentials of $f _k$ with respect the $\natural$-metrics satisfies
\begin{equation}
\sum_k||\wedge^2{\rm d}f_k||< 1.\label{eq-a}
\end{equation}
(Notice that $\sum_k||\wedge^2{\rm d}f_k||= 1$ if $X=Y$ and $f$ is the identity map.)

If either the map $f$ is {\it spin} or the {\it universal covering of $Y$ is spin}, then the {\it topological degree of $f$ is zero}.

{\it Proof.} First, let $X$ and $Y$ be spin, let $f\colon Y\to X$ be a smooth map and
 let $V\to Y$ be the $f$-pullback of the spin bundle $\mathcal S (X)\to X$ to $Y$.
Then, if $\mathcal R^X\geq 0$ and $ f\colon Y\to X$ is $\natural$-area decreasing at point $y\in Y$, i.e., $\|{\wedge}^2{\rm d}f(y)\|\leq 1$, then
according to~\cite{GS2002} (also see~\cite{Li2010}) the lowest eigenvalue of the operator
$\mathcal R_{\otimes V} $ at $y\in Y$ satisfies
\[
\lambda_{\otimes V}\geq -{\operatorname{Sc}(Y,y)\over 4},
\]
where this inequality is strict if $ f$ is strictly $\natural$-area decreasing at $y$.

 Next, let $X_k$ be spin, let $X=\bigtimes_kX_k$ and let $V\to Y$ be the tensor product $V=\bigotimes_kV_k$
of the pullbacks $V_k=f^\ast_k (\mathbb S_k)\to Y$ of $f^\ast_k (\mathbb S_k)=\mathbb S(X_k)\to X_k$ to $Y$ for $f_k\colon Y\to X_k$.

Then, if the maps $f_k$ are $\natural$-area decreasing and at least
one of $f_k$ is strictly $\natural$-area decreasing at a point $y\in Y$
 then, assuming $Y$ is connected and spin, the Dirac operator $\mathcal D_{\otimes V}$ on $Y$ has {\it index zero}.

On the other hand, if $\chi(X_k)\neq 0$, if $\dim(Y)=\dim(X)$ and
$\deg (f)\neq 0$,\footnote{If $\dim(Y)=\dim(X)+4m$, then instead of $\deg(f)\neq 0$ one assumes that the $f$-pullback of a generic point $x\in X$ has non-zero
hat A-genus.}
 then $\operatorname{ind}( \mathcal D_{\otimes V})\neq 0$ by the
 Atiyah--Singer theorem (compare with \cite{GS2002, Li2010, Ll1998}).

This proves \S5.A in the case where the manifold $X$ is spin and it contains neither a $Z$-factor, nor an
odd spherical factor.

 To pass to the general case we argue as follows:

1. Odd dimensional spheres are suspended to even dimensional
ones $S^{n_k}\leadsto S^{n_k+1}$, where these suspensions
 are accompanied by multiplying $Y$ by a long circles and a suspending $[f_k\colon Y\to S^{n_k}\leadsto [Y\times S^1\to S^{n_k+1}]$ as
 in \cite{Ll1998}, also see \cite[\S3.4.1]{Gr2021} and \cite{Gr2023}.

2. If a $Z$, which may be assumed even-dimensional, is enlargeable,
it supports an almost flat bundle, say $W \to Z$ with non-zero top-dimensional Chern class and the above $V\to Y$ is tensored by the
pullback $f_Z^\ast(W)\to Y$ of $W$ to $Y$.

3. If neither $X$ nor $Y$ are spin but the map $f$ is spin, then
the Dirac operator $\mathcal D_{\otimes V}$ is defined (this is explained in the present context in~\cite{Mi1995} and in \cite{GS2002}) and the above applies.

\textbf{5.B. Spherical trace and symplectic remarks.} The $ ||\wedge^2{\rm d}f_k||$ contribution of each spherical factor $X_k$ with constant sectional curvature can be replaced
 in the formula~\eqref{eq-a} by an a priori smaller
entity, that is, $ 2||\wedge^2{\rm d}f||_{\rm trace}\over n_k(n_k-1)$,
where
\[
 ||\wedge^2{\rm d}f_k(y)||_{\rm trace}=\sum_{1\leq i<j\leq n+N} \lambda_{i,k}(y) \lambda_{j,k}(y),
 \]
and where the numbers $\lambda_{j,k}(y)\geq 0$ are defined by diagonalizing the
differential ${\rm d}f_k\colon T_y(Y)\to T_{f_k}(X_k)$
with an {\it orthonormal} frame $e_{i,k}\in T(y)(Y)$,
which is sent by ${\rm d}f_k$ to an {\it orthogonal} frame in $T_{f_k}(X_k)$ with the vectors of
lengths $\lambda_{j,k}(y)$.

In fact, this follows from~\cite[formula~(4.6)]{Ll1998} (also \cite{Li2010} and~\cite[\S3.4]{Gr2021}).

The $S^2$ factors in $X$ contribute to complex line bundles as
$\otimes$-factors in $ V\to Y$.
This, in view of Schr\"odinger--Hitchin (see~\cite{Hi1974}) formula for $\mathcal D_{\otimes L}$ allows one to replace the product of these
$S^2$ by a single (quasi)symplectic manifold (compare with \cite[\S2.7 and \S3.4.4(4)]{Gr2021}).

\textbf{5.C. $\boldsymbol{\operatorname{Sc}_{{\rm area},\widetilde {\rm sp}}^{\rtimes_\downarrow}}$-additivity corollary.}
Let $X_k$ be manifolds as in $5.A$, where we additionally assume that they are spin and have constant scalar curvatures.
Then the fundamental classes~$[X_k]$ satisfy the $\widetilde {\rm sp}$-version of the
$\operatorname{Sc}_{{\rm sup.area}}^{\rtimes_\downarrow}$-additivity \eqref{eq-c} in \S4.F:
\[
\operatorname{Sc}_{{\rm sup.area}, \widetilde{\rm sp}}^{\rtimes_\downarrow}( \bigtimes_k X_k)=
\sum_k\operatorname{Sc}_{ \widetilde{\rm sp}}^{\rtimes_\downarrow}(X_k)=
\sum_k\operatorname{Sc}^{\rtimes}(X_k)=\sum_k\operatorname{Sc}(X_k).
\]
Consequently,
\[
\operatorname{Sc}_{\widetilde{\rm sp}}^{\rtimes_\downarrow}( \bigtimes_{\sup_k} X_k)
=\sum_k\operatorname{Sc}(X_k).
\]

 \textbf {5.D. Questions.} (i) Does vanishing of $\operatorname{Sc}^{\rtimes_\downarrow} [X_k]_{\mathbb Q}$ for closed manifolds $X_k$ (this is a homotopy invariant of $X$)
 implies vanishing of $\operatorname{Sc}^{\rtimes_\downarrow} [\bigtimes_kX_k]_{\mathbb Q}$?

There are examples of manifolds $X_k$, where
$\operatorname{Sc}^{\rtimes_\downarrow} [X_i]=0$ and where their products admit metrics with $\operatorname{Sc}>0$; hence,
 $\operatorname{Sc}^{\rtimes_\downarrow} [\bigtimes_iX_i]>0 $ for these $X_i$, see~\cite{GH2023}.

(ii) Do products of
spheres $X=S^{n_1}\times S^{n_2}$, $n_1,n_2\geq 2$, admit Riemannian metrics $g_\varepsilon$, for all $\varepsilon>0$, with $\operatorname{Sc}(g_\varepsilon)\geq 1$ and such that {\it all non-zero} homology classes $h$ in $H_{n_1}(X)$ and in $H_{n_2}(X)$ satisfy
$\operatorname{Sc}_{\operatorname{area}_{ g_\varepsilon}}^{\rtimes_\downarrow} (h)\leq \varepsilon$?

The existence of such a $g_\varepsilon$, for $n_1=n_2=2$, would imply
the absence of the lower bounds on the 2-{\it systoles} of manifolds $(X, g)$ in terms of $\sigma(g)=\inf_{x\in X} \operatorname{Sc} (X,g,x)>0$,\footnote{Such a counter example would undermine (but not disprove) the conjectural bound ${\rm waist}_2(X)\leq \frac{{\rm const}_n}{\sigma}$ for compact Riemannian $n$-manifolds with $\operatorname{Sc}(X)\geq \sigma>0$. Thus, it may be safer
to assume $n_1,n_2\geq 3$.}
\[
\sup_{\sigma(g)\geq 1} {\rm syst}_2(X, g)= \infty.
\]

Recall, that the 2-{\it systole} is the {\it infimum of the areas of all non-zero classes} $h\in H_2(X)$, for $\operatorname{area}(h)=\inf_{[c]\in h} \operatorname{area}(c)$ for the 2-cycles $c\subset Y$ that represent~$h$.\footnote{There are bounds on the 2-systoles of manifolds $X$ with $\operatorname{Sc}^\rtimes (X)\geq \sigma$ in terms of their $\tilde \square^\perp$-spreads (see
\cite{Ri2020, Ze2020}) which are proved as $\square^{\exists\exists}(n,m)$-inequality in~\S2.B with a use of minimal hypersurfaces and $\mu$-bubbles.
 Also there are similar bounds on the {\it stable systoles} of spin manifolds obtained with Dirac operators twisted with line bundles, where, recall,
 ${\rm st.syst}_2(X)=\liminf_{N\to \infty} \frac{\operatorname{area}(Nh)}{N}$.}

 (iii) Let $X$ be a compact symmetric space. What is the {\it minimal}
 seminorm on the linear maps $\lambda^2{\rm d}\colon  \mathbb \wedge^2R^n\to \wedge^2T(X)$,
say $||\lambda^2{\rm d}||_{\min}$ such that
 the $\natural$-normalized inequality $||\wedge^2{\rm d}f||_{\min} < 1$ for smooth equidimensional spin maps $f\colon Y\to X$ would imply
 that $\deg(f)=0$?\footnote{This norm must be invariant under isometries of $X$.}

 (If $X$ is the products of spheres, this seminorm is equal to the sum of the mean trace norms (as in \S5.B) for maps
 $\mathbb R^n\overset {{\rm d}_k}\to X_k=S^{n_k}$ and for
 all symmetric spaces $X$ of dimension $\geq 4$ with $\chi(X)\neq 0$ this norm is, probably, strictly smaller then the sup-norm $||\wedge^2{\rm d}||$ from the Goette--Semmelmann theorem.)

\section[P-families of maps to product of spheres]{$\boldsymbol{P}$-families of maps to product of spheres}\label{section6}

Let $Y=(Y,g) $ be an $n$-dimensional Riemannian manifold with $\operatorname{Sc}(X)>0$, where as earlier
$g^\natural =\operatorname{Sc}(Y)\cdot g$,
 let $h_m\in H_m(Y)$ be a homology class and let $P$ be a locally contractible topological space, e.g., a manifold and $h_K\in H_K(P)$ be a homology class.

Let $X$ be a product of spheres of variable radii,
\[
X=\bigtimes _kS^{n_k}(R_k),
\]
where $\dim(X)=\sum_kn_k=m+K$, and where the spheres are endowed with the usual metrics with sectional curvatures $1/R_k^2$.

 Let $F\colon Y \times P\to X$ be a continuous map such that
 the maps $F_p=F_{|Y\times p}\colon Y\to X$ are smooth and $C^1$-continuous in $p\in P$.

Let the universal covering of $Y$ be spin and let $h_m$ be equal to the homology class of the pullback of a genetic point under a smooth map $\phi\colon Y\to Z$, where $Z$ is a smooth enlargeable manifold of
 dimension $\dim(Y)-m$. For instance, $m=n$ and $h_m=[Y]$ or
 $Y =Y_0^m\times \mathbb T^{n-m}$ and $h_m=[Y_0]$.

\textbf {6.A. Theorem.}  Let the norms of the exterior squares of the differentials of the maps $f_k\colon Y\to S^{n_k}(R_k)$ with respect to the $\natural$-metrics in $Y$ and in $S^{n_k}(R_k)$ satisfy
\begin{equation*}
\sum_k||\wedge^2{\rm d}f_k||< 1. 
\end{equation*}
 Then, in the following two cases,
the $F$-image $F_\ast(h_m\otimes h_K)\in H_{m-K}(X)=\mathbb Z$ vanishes:
\begin{enumerate}\itemsep=0pt
\item[(1)] The ranks of the (differentials of the) maps $f_{p}\colon Y\to X$ are everywhere $\leq m$, e.g., $\dim(Y)=m$ and $h_m=[Y]$.
\item[(2)] The dimension of $Y$ is bounded by $n\leq 8$.
\end{enumerate}

 {\it Proof.} Case~1. If $\operatorname{ranks}(f_p)\leq m$, then the Llarull (Listing) trace inequality~(4.6) in \cite{Ll1998} together withe the above $\lambda_{\otimes V}$-additivity show that index of the family of the Dirac operators on~$Y$, twisted with the pullbacks of $\bigotimes_k \mathbb S_k$ as in \S5.A, vanishes and the Atiyah--Singer theorem for families shows that $F_\ast(h_m\otimes h_K)=0$. (See in \cite[\S4]{Gr2021} and references therein.)

 Case~2. If $\dim(Y)\leq 8$, then, at last generically, the homology class
 $h_m$ can be realized by an $m$-submanifold $Y_0\subset Y $ such that the product $Y_0\times \mathbb T^{\dim(Y)-m}$ admits a warped product metric~$g^\rtimes $ such that $\operatorname{Sc} (g^\rtimes, y)\geq \operatorname{Sc} (Y,y)$
 for all $y\in Y_0$ (see \cite[\S3]{Gr2023} and references therein).
 Now the case~1 applies to $Y_0\times \mathbb T^{\dim(Y)-m}$ and the proof follows.

{\it Remarks/Problems.} (a) For all we know, the spin and $\dim(Y)\leq 8$ condition are redundant and there is a fair chance that a further study of singularities of minimal hypersurfaces in he spirit of~\cite{SY2017} and/or~\cite{Lo2018,Lo2018a} will allow one to remove the latter. But removing the spin condition needs a new idea.

The argument in Case~1 can be extended to maps of foliated manifolds to $\bigtimes_k S^{n_k}$ as in~\cite{SWZ2021}, but a foliated version of Case~2 is problematic.

\section{Spinorial curvature}\label{section7}

Given a closed orientable even dimensional Riemannian manifold $Y$ let
$\mathbb S{\sf p}.\mathrm{curv}(Y)$ be the infimum of the numbers $\kappa\geq 0$ such that
there exist
a complex vector bundle $V\to X$ with a~unitary connection such that
\[
\bigl(\mathrm{Ch}(V)\smile \hat A\bigr)[Y]\neq 0
\]
 and
 the lowest eigenvalue of the operators $\mathcal K_{\otimes V}$ on $(\mathcal S\otimes V)_y$ (see Section~\ref{section5}) satisfies
\[
\lambda_{\otimes V}\geq -\kappa
\]
 at all points $y\in Y$.\footnote{Although the spin bundle $\mathcal S\to Y$ is defined only for spin manifolds, this
definition, being local, makes sense for all $Y$, since $\lambda_{\otimes V}$ does not depend on the spin structure.}

Observe that
\[
\mathbb S\mbox{\sf p}.{\rm curv}(Y_1\times Y_2)\leq \mathbb S\mbox{\sf p}.{\rm curv}(Y_1)+ \mathbb S\mbox{\sf p}.{\rm curv}( Y_2)
\]
 by
the inequality $\lambda_{\otimes 1\otimes 2}\geq \lambda_{\otimes 1}+\lambda_{\otimes 2} $
from Section~\ref{section5},
 that
\[
\mathbb S\mbox{\sf p}.{\rm curv}(Y)=0
\]
for enlargeable manifolds $Y$
and that if the universal coverings of $Y$ is {\it spin}, then
\[
\operatorname{Sc}^\rtimes(Y)\leq 4\kappa
\]
by the ($\mathbb T^\rtimes$-stabilized) index theorem, SLWB-formula and Kato's inequality.\footnote{Instead of the $\mathbb T^\rtimes$-stabilization and Kato's inequality one may use Kazdan--Warner conformal change theorem~\cite{KW1975} and conformal invariance of harmonic spinors \cite{Hi1974}.}

{\it Remarks on $\lambda_1(X, \beta)$ for $\beta<1/4$.} (a) The refined Kato inequality strengthens the above to
\[
\lambda_1(X,\beta)\leq 4\kappa \qquad \text{for}\quad \beta= (n-1)/4n,
\]
where $\lambda_1(X,\beta)$ is the lowest eigenvalue of $-\Delta+\beta \operatorname{Sc}(X)$ (see~\S1.D).

(b) The Kazdan--Warner conformal change theorem~\cite{KW1975}
and conformal invariance of harmonic spinors~\cite{Hi1974} show that
 if $\lambda_1(X,\beta)>0$ for $\beta=(n-2)(4(n-2)$, then
$X$ supports no non-zero harmonic spinors.

However, it is unclear how to extract further geometric, rather than topological information from
 the inequality $\lambda_1(X,\beta)> \sigma $ for $\beta < (n-1)/4n$ and
 $\sigma>0$.

Next, let $X$ be a Riemannian manifold, let $h_m\in H_m(X)$
be a homology class and let $\mathcal Y $ be a~class of smooth closed orientable $m$-manifolds $Y$ along with maps
$f\colon Y\to X$.

 Define $ \mathbb S${\sf p}$.{\rm curv}^\downarrow _{\mathcal Y}(X)$ via smooth maps
$F \colon Y\times \mathbb T^N\to X\times \mathbb T^N$ and Riemannian metrics $G$ on $Y\times \mathbb T^N$ as the infimum
\[
\mathbb S\mbox{\sf p}.{\rm curv}^\downarrow _{\mathcal Y}(X)\mbox{\sf p}\inf_{ Y,N,G,F} \mathbb S\mbox{\sf p}.{\rm curv} \bigl(Y\times \mathbb T^N,G\bigr),
\]
where the infimum is taken over $N$ such that $m+N$ is {\it even}, where $F$ is {\it area decreasing} with respect to the metric $G$ and where
\[
F_\ast\bigl[Y\times\mathbb T^N\bigr]=h_m\otimes \bigl[\mathbb T^N\bigr]\in H_{m+N}\bigl(X\times \mathbb T^N\bigr) \qquad \text{and} \qquad (Y, F_{|Y\times 0})\in \mathcal Y.
\]

 Clearly, by the above,
if the universal coverings of manifolds $Y\in \cal Y$ are spin,\footnote{This condition is necessary but its $\mathbb Q$-version may be redundant.} then
\[
\mathbb S\mbox{\sf p}.{\rm curv}^\downarrow _{\mathcal Y}(h_m)\geq\frac{1}{4} \bigl(\operatorname{Sc}^{\rtimes _\downarrow}_{{\rm area},\mathcal Y }(h_m)\bigr).
\]

{\it Remark.} If the universal coverings of the manifolds $Y\in \cal Y$ are spin, then the fundamental classes $[X]$ of compact symmetric spaces $X$ with $\chi(X)\neq 0$, satisfy the equally
\[
\mathbb S\mbox{\sf p}.{\rm curv}^\downarrow _{\mathcal Y}(X)=\frac{1}{4} \bigl(\operatorname{Sc}^{\rtimes _\downarrow}_{{\rm area},\mathcal Y }(h)\bigr)
\]
 by the $\mathbb T^\rtimes$-stabilized Goette--Semmelmann theorem and this equally applies to products $X=\bigtimes X_i$, where $X_i$ are as in \S5.A.

Possibly, (a version of) this equality holds true for all symmetric spaces but it seems unlikely in general, even for rational homology classes $h$, that the Dirac operator
 is the only source of bounds on $\operatorname{Sc}^{\rtimes _\downarrow}_{{\rm area},\mathcal Y }(h).$

{\it Power stabilization.} Let
\begin{gather*}
X^M= \underset {M}{ \underbrace {X\times X\times \cdots \times X}},\\
\operatorname{Sc}^{\rtimes _\downarrow}_{{\rm area},\mathcal Y }\bigl(\otimes ^{[\infty/\infty]}h_m\bigr)=\sup_{M=1,2,\dots}\frac {1}{M} \operatorname{Sc}^{\rtimes _\downarrow}_{{\rm area},\mathcal Y }(h_{Mm}),\qquad  h_{Mm}=
\otimes^M h_m\in H_{Mm}\bigl(X^M\bigr)
\end{gather*}
and
\begin{gather*}
\mathbb S\mbox{\sf c}.{\rm curv}^{\rtimes _\downarrow}_{{\rm area},\mathcal Y }\bigl(\otimes ^{[\infty/\infty]}h_m\bigr)=\inf_{M=1,2,\dots }\frac {1}{M} \mathbb S\mbox{\sf c}.{\rm curv}^{\rtimes _\downarrow}_{{\rm area},\mathcal Y }(h_{Mm}),\\
h_{Mm}= \otimes^M h_m\in H_{Mm}\bigl(X^M\bigr).
\end{gather*}

{\it Questions.} \textbf{I.} What are further instances (besides the above $h=[X]$) of the equality
\[
\mathbb S\mbox{\sf p}.{\rm curv}^\downarrow _{\mathcal Y}(\otimes ^{[\infty/\infty]}h_m)=\frac{1}{4} \bigl(\operatorname{Sc}^{\rtimes _\downarrow}_{{\rm area},\mathcal Y }\bigl(\otimes ^{[\infty/\infty]}h_m\bigr)\bigr),
\]
and what are examples where this fails to be true?

 \textbf{II.} Can one pass to the limit, set $M=\infty$ and prove scalar curvature bounds for
``Riemannian metrics'' $G$ on infinite dimensional manifolds $X$, e.g., where such a G differs from the infinite sum of Riemannian metrics, $\sum_1^\infty g_i$, on $X=\bigtimes_1^\infty(X_i,g_i)$ (and/or on $Y$ mapped to $X$) by a ``fast decaying in
$i$'' error term $\Delta$?

 {\it Remarks.} (a) If $\Delta=\Delta_{i,j}$ decays very fast, for $i$ and/or $j$ tending to infinity, then finite products $X_M=\bigtimes_1^M X_i$ embed to
$X= \bigtimes_1^\infty(X_i,g_i)$ with small relative curvatures and a bound on ``$\operatorname{Sc}(X)$'' may be derived in some cases from such a bound on $X_M$, but it would be more interesting to develop a
 truly infinite dimensional argument for bounds on ``$\operatorname{Sc}(X)$'' and/or to find applications of such
 bounds.

 {\it Test question.} Let $X=\{x_i\}_{\sum_i x_i^2\leq\infty}$ be the Hilbert space and $G=G_{ij}$ be a smooth Riemannian metric on $X$, which is {\it greater} than the background Hilbertian metric,{\samepage
\[
G(\tau,\tau)\geq ||\tau||^2
\]
 for all tangent vectors $ \tau\in T(X)$ and let $M=3,4,\dots$ be an integer.}

 Can the $M$-{\it scalar curvature} of $G$ (defined below) {\it be strictly positive}, say
$\operatorname{Sc}_M(G)\geq 1$?
 Here $\operatorname{Sc}_M$ is the function on the tangent $M$-planes
 $P^M\in T_x(X)$, $x\in X$, which is equal to the scalar curvature at zero in $P^M$ $\bigl({=}\,\mathbb R^M\bigr)$ of the Riemannian metric induced by the exponential map $\exp\colon P^M\to X$ from~$G$. (It may be worthwhile to compare $\operatorname{Sc}_M$ with
 with the {\it $m$-intermediate curvature} from~\cite{BHG2022}.)

(b) A natural approach to these problems is by a finite-dimensional approximation as in (a) but this seems that uncomfortably restrictive conditions on $G$
are needed (compare with~\cite{Gr2014}).

(c) Basic features of positive scalar curvature have their counterparts
for {\it mean convex} hypersurfaces (see~\cite{Gr2019}), where the infinite dimensional geometry is a bit more transparent than that of the scalar curvature.

\subsection*{Acknowledgements}
I am thankful to anonymous referees for helpful suggestions and for pointing out errors in the first version of this paper.

\pdfbookmark[1]{References}{ref}
\LastPageEnding


\begin{thebibliography}{99}
\footnotesize\itemsep=0pt

\bibitem{BF2017}
Borisov D., Freitas P., The spectrum of geodesic balls on spherically symmetric
  manifolds, \href{https://doi.org/10.4310/CAG.2017.v25.n3.a1}{\textit{Comm. Anal. Geom.}} \textbf{25} (2017), 507--544,
  \href{https://arxiv.org/abs/1603.02399}{arXiv:1603.02399}.

\bibitem{BR2021}
Botvinnik B., Rosenberg J., Positive scalar curvature on manifolds with fibered
  singularities, \href{https://doi.org/10.1515/crelle-2023-0055}{\textit{J.~Reine Angew. Math.}} \textbf{803} (2023), 103--136,
  \href{https://arxiv.org/abs/1808.06007}{arXiv:1808.06007}.

\bibitem{BHG2022}
Brendle S., Hirsch S., Johne F., A generalization of {G}eroch's conjecture,
  \href{https://doi.org/10.1002/cpa.22137}{\textit{Comm. Pure Appl. Math.}} \textbf{77} (2024), 441--456,
  \href{https://arxiv.org/abs/2207.08617}{arXiv:2207.08617}.

\bibitem{BH2009}
Brunnbauer M., Hanke B., Large and small group homology, \href{https://doi.org/10.1112/jtopol/jtq014}{\textit{J.~Topol.}}
  \textbf{3} (2010), 463--486, \href{https://arxiv.org/abs/0902.0869}{arXiv:0902.0869}.

\bibitem{CHS2023}
Cecchini S., Hanke B., Schick T., Lipschitz rigidity for scalar curvature,
  \href{https://arxiv.org/abs/2206.11796}{arXiv:2206.11796}.

\bibitem{CZ2021a}
Cecchini S., Zeidler R., Scalar and mean curvature comparison via the {D}irac
  operator, \textit{Geom. Topol.}, {t}o appear, \href{https://arxiv.org/abs/2103.06833}{arXiv:2103.06833}.

\bibitem{CF1974}
Chavel I., Feldman E.A., The first eigenvalue of the {L}aplacian on manifolds
  of non-negative curvature, \textit{Compositio Math.} \textbf{29} (1974),
  43--53.

\bibitem{Ch1979}
Cheeger J., On the spectral geometry of spaces with cone-like singularities,
  \href{https://doi.org/10.1073/pnas.76.5.2103}{\textit{Proc. Nat. Acad. Sci. USA}} \textbf{76} (1979), 2103--2106.

\bibitem{CMS2023}
Chodosh O., Mantoulidis C., Schulze F., Generic regularity for minimizing
  hypersurfaces in dimensions~9 and~10, \href{https://arxiv.org/abs/2302.02253}{arXiv:2302.02253}.

\bibitem{F-CS1980}
Fischer-Colbrie D., Schoen R., The structure of complete stable minimal
  surfaces in {$3$}-manifolds of non-negative scalar curvature, \href{https://doi.org/10.1002/cpa.3160330206}{\textit{Comm.
  Pure Appl. Math.}} \textbf{33} (1980), 199--211.

\bibitem{GS1999}
Goette S., Semmelmann U., {${\rm Spin}^c$} structures and scalar curvature
  estimates, \href{https://doi.org/10.1023/A:1013035721335}{\textit{Ann. Global Anal. Geom.}} \textbf{20} (2001), 301--324,
  \href{https://arxiv.org/abs/math.DG/9905089}{arXiv:math.DG/9905089}.

\bibitem{GS2002}
Goette S., Semmelmann U., Scalar curvature estimates for compact symmetric
  spaces, \href{https://doi.org/10.1016/S0926-2245(01)00068-7}{\textit{Differential Geom. Appl.}} \textbf{16} (2002), 65--78,
  \href{https://arxiv.org/abs/math.DG/0010199}{arXiv:math.DG/0010199}.


\bibitem{Gr1996}
Gromov M., Positive curvature, macroscopic dimension, spectral gaps and higher
  signatures, in Functional {A}nalysis on the {E}ve of the 21st {C}entury,
  {V}ol.~{II} ({N}ew {B}runswick, {NJ}, 1993), \textit{Progr. Math.}, Vol.~132,
  \href{https://doi.org/10.1007/978-1-4612-4098-3_1}{Birkh\"auser}, Boston, MA, 1996, 1--213.

\bibitem{Gr2014}
Gromov M., Manifolds: {W}here do we come from? {W}hat are we? {W}here are we
  going?, in The {P}oincar\'e {C}onjecture, \textit{Clay Math. Proc.}, Vol.~19,
  American Mathematical Society, Providence, RI, 2014, 81--144.

\bibitem{Gr2018'}
Gromov M., Hilbert volume in metric spaces, \href{https://arxiv.org/abs/1811.04332}{arXiv:1811.04332}.

\bibitem{Gr2019}
Gromov M., Mean curvature in the light of scalar curvature, \href{https://doi.org/10.5802/aif.3347}{\textit{Ann. Inst.
  Fourier}} \textbf{69} (2019), 3169--3194, \href{https://arxiv.org/abs/1812.09731}{arXiv:1812.09731}.

\bibitem{Gr2021}
Gromov M., Four lectures on scalar curvature, in Perspectives in {S}calar
  {C}urvature. {V}ol.~1, \href{https://doi.org/10.1142/9789811273223_0001}{World Scientific Publishing}, Hackensack, NJ, 2023,
  1--514, \href{https://arxiv.org/abs/1902.10612}{arXiv:1902.10612}.

\bibitem{Gr2023}
Gromov M., Scalar curvature, injectivity radius and immersions with small
  second fundamental, \href{https://arxiv.org/abs/2203.14013}{arXiv:2203.14013}.

\bibitem{GH2023}
Gromov M., Hanke B., Torsion obstructions to positive scalar curvature,
  \href{https://arxiv.org/abs/2112.04825}{arXiv:2112.04825}.

\bibitem{GL1980}
Gromov M., Lawson Jr. H.B., Spin and scalar curvature in the presence of a
  fundamental group.~{I}, \href{https://doi.org/10.2307/1971198}{\textit{Ann. of Math.}} \textbf{111} (1980), 209--230.

\bibitem{GL1983}
Gromov M., Lawson Jr. H.B., Positive scalar curvature and the {D}irac operator
  on complete {R}iemannian manifolds, \href{https://doi.org/10.1007/BF02953774}{\textit{Inst. Hautes \'Etudes Sci. Publ.
  Math.}} \textbf{58} (1983), 83--196.

\bibitem{GZ2021}
Gromov M., Zhu J., Area and {G}auss--{B}onnet inequalities with scalar
  curvature, \href{https://doi.org/10.4171/cmh/570}{\textit{Comment. Math. Helv.}} \textbf{99} (2024), 355--395,
  \href{https://arxiv.org/abs/2112.07245}{arXiv:2112.07245}.

\bibitem{HKRS}
Hanke B., Kotschick D., Roe J., Schick T., Coarse topology, enlargeability, and
  essentialness, \href{https://doi.org/10.24033/asens.2073}{\textit{Ann. Sci. \'Ec. Norm. Sup\'er.~(4)}} \textbf{41}
  (2008), 471--493, \href{https://arxiv.org/abs/0707.1999}{arXiv:0707.1999}.

\bibitem{He2000}
Herzlich M., Refined {K}ato inequalities in {R}iemannian geometry, in
  Journ\'ees ``\'Equations aux {D}\'eriv\'ees {P}artielles'' ({L}a {C}hapelle
  sur {E}rdre, 2000), Universit\'e de Nantes, Nantes, 2000, 1--11.

\bibitem{Hi1974}
Hitchin N., Harmonic spinors, \href{https://doi.org/10.1016/0001-8708(74)90021-8}{\textit{Adv. Math.}} \textbf{14} (1974), 1--55.

\bibitem{KW1975}
Kazdan J.L., Warner F.W., Scalar curvature and conformal deformation of
  {R}iemannian structure, \href{https://doi.org/10.4310/jdg/1214432678}{\textit{J.~Differential Geometry}} \textbf{10} (1975),
  113--134.

\bibitem{LM1989}
Lawson Jr. H.B., Michelsohn M.-L., Spin geometry, \textit{Princeton Math. Ser.},
  Vol.~38, \href{https://doi.org/10.1515/9781400883912}{Princeton University Press}, Princeton, NJ, 1989.

\bibitem{LM2021}
Liokumovich Y., Maximo D., Waist inequality for 3-manifolds with positive
  scalar curvature, in Perspectives in {S}calar {C}urvature. {V}ol.~2, \href{https://doi.org/10.1142/9789811273230_0022}{World
  Scientific Publishing}, Hackensack, NJ, 2023, 799--831, \href{https://arxiv.org/abs/2012.12478}{arXiv:2012.12478}.

\bibitem{Li2010}
Listing M., Scalar curvature on compact symmetric spaces, \href{https://arxiv.org/abs/1007.1832}{arXiv:1007.1832}.

\bibitem{Ll1998}
Llarull M., Sharp estimates and the {D}irac operator, \href{https://doi.org/10.1007/s002080050136}{\textit{Math. Ann.}}
  \textbf{310} (1998), 55--71.

\bibitem{Lo2018}
Lohkamp J., Minimal smoothings of area minimizing cones, \href{https://arxiv.org/abs/1810.03157}{arXiv:1810.03157}.

\bibitem{Lo2018a}
Lohkamp J., Contracting maps and scalar curvature, \href{https://arxiv.org/abs/1812.11839}{arXiv:1812.11839}.

\bibitem{MN2011}
Marques F.C., Neves A., Rigidity of min-max minimal spheres in three-manifolds,
  \href{https://doi.org/10.1215/00127094-1813410}{\textit{Duke Math.~J.}} \textbf{161} (2012), 2725--2752, \href{https://arxiv.org/abs/1105.4632}{arXiv:1105.4632}.

\bibitem{Mi1995}
Min-Oo M., Scalar curvature rigidity of certain symmetric spaces, in Geometry,
  {T}opology, and {D}ynamics ({M}ontreal, {PQ}, 1995), \textit{CRM Proc.
  Lecture Notes}, Vol.~15, \href{https://doi.org/10.1090/crmp/015/08}{American Mathematical Society}, Providence, RI, 1998,
  127--136.

\bibitem{O'N1966}
O'Neill B., The fundamental equations of a submersion, \href{https://doi.org/10.1307/mmj/1028999604}{\textit{Michigan
  Math.~J.}} \textbf{13} (1966), 459--469.

\bibitem{QW1999}
Qu C.K., Wong R., ``{B}est possible'' upper and lower bounds for the zeros of
  the {B}essel function {$J_\nu(x)$}, \href{https://doi.org/10.1090/S0002-9947-99-02165-0}{\textit{Trans. Amer. Math. Soc.}}
  \textbf{351} (1999), 2833--2859.

\bibitem{Ri2020}
Richard T., On the 2-systole of stretched enough positive scalar curvature
  metrics on {$\mathbb S^2\times\mathbb S^2$}, \href{https://doi.org/10.3842/SIGMA.2020.136}{\textit{SIGMA}} \textbf{16}
  (2020), 136, 7~pages, \href{https://arxiv.org/abs/2007.02705}{arXiv:2007.02705}.

\bibitem{SY1979}
Schoen R., Yau S.-T., Existence of incompressible minimal surfaces and the
  topology of three-dimensional manifolds with non-negative scalar curvature,
  \href{https://doi.org/10.2307/1971247}{\textit{Ann. of Math.}} \textbf{110} (1979), 127--142.

\bibitem{SY2017}
Schoen R., Yau S.-T., Positive scalar curvature and minimal hypersurface
  singularities, in {D}ifferential {G}eometry, {C}alabi--{Y}au {T}heory, and
  {G}eneral {R}elativity. {P}art~2, \textit{Surv. Differ. Geom.}, Vol.~24,
  \href{https://doi.org/10.4310/SDG.2019.v24.n1.a10}{International Press}, Boston, MA, 2022, 441--480, \href{https://arxiv.org/abs/1704.05490}{arXiv:1704.05490}.

\bibitem{Sm2003}
Smale N., Generic regularity of homologically area minimizing hypersurfaces in
  eight-dimensional manifolds, \href{https://doi.org/10.4310/CAG.1993.v1.n2.a2}{\textit{Comm. Anal. Geom.}} \textbf{1} (1993),
  217--228.

\bibitem{St1992}
Stolz S., Simply connected manifolds of positive scalar curvature, \href{https://doi.org/10.2307/2946598}{\textit{Ann.
  of Math.}} \textbf{136} (1992), 511--540.

\bibitem{SWZ2021}
Su G., Wang X., Zhang W., Nonnegative scalar curvature and area decreasing maps
  on complete foliated manifolds, \href{https://doi.org/10.1515/crelle-2022-0038}{\textit{J.~Reine Angew. Math.}} \textbf{790}
  (2022), 85--113, \href{https://arxiv.org/abs/2104.03472}{arXiv:2104.03472}.

\bibitem{WXY2021}
Wang J., Xie Z., Yu G., A proof of {G}romov's cube inequality on scalar
  curvature, \href{https://arxiv.org/abs/2105.12054}{arXiv:2105.12054}.

\bibitem{Ze2020}
Zeidler R., Width, largeness and index theory, \href{https://doi.org/10.3842/SIGMA.2020.127}{\textit{SIGMA}} \textbf{16}
  (2020), 127, 15~pages, \href{https://arxiv.org/abs/2002.13754}{arXiv:2002.13754}.

\bibitem{Ze2019}
Zeidler R., Band width estimates via the {D}irac operator,
  \href{https://doi.org/10.4310/jdg/1668186790}{\textit{J.~Differential Geom.}} \textbf{122} (2022), 155--183,
  \href{https://arxiv.org/abs/1905.08520}{arXiv:1905.08520}.

\bibitem{Zhu2019}
Zhu J., Rigidity of area-minimizing {$2$}-spheres in {$n$}-manifolds with
  positive scalar curvature, \href{https://doi.org/10.1090/proc/15033}{\textit{Proc. Amer. Math. Soc.}} \textbf{148}
  (2020), 3479--3489, \href{https://arxiv.org/abs/1903.05785}{arXiv:1903.05785}.

\bibitem{Zhu2020}
Zhu J., Rigidity results for complete manifolds with nonnegative scalar
  curvature, \href{https://doi.org/10.4310/jdg/1701804153}{\textit{J.~Differential Geom.}} \textbf{125} (2023), 623--644,
  \href{https://arxiv.org/abs/2008.07028}{arXiv:2008.07028}.

\end{thebibliography}
\end{document}